\RequirePackage{ifpdf}
\ifpdf 
\documentclass[pdftex]{sigma}
\else
\documentclass{sigma}
\fi

\usepackage{array}
\usepackage{euscript}
\usepackage{euscript}
\usepackage{fontenc}
\usepackage[all]{xy}
\usepackage{xspace}

\numberwithin{equation}{section}

\newtheorem{Theorem}{Theorem}[section]
\newtheorem{Lemma}[Theorem]{Lemma}
\newtheorem{Corollary}[Theorem]{Corollary}
\newtheorem{Conjecture}[Theorem]{Conjecture}
\newtheorem{Proposition}[Theorem]{Proposition}
{\theoremstyle{definition}
 \newtheorem{Example}[Theorem]{Example}
\newtheorem{Definition}[Theorem]{Definition}
\newtheorem{Remark}[Theorem]{Remark}
}

\newcommand{\C}[1]{\ensuremath{\mathbb{C}^{#1}}}
\newcommand{\R}[1]{\ensuremath{\mathbb{R}^{#1}}}
\newcommand{\Z}[1]{\ensuremath{\mathbb{Z}^{#1}}}

\newcommand{\OO}[1]{
  \ensuremath{
    \mathcal{O}
    \ifthenelse{\equal{#1}{0}}
      {}
      {\left({#1}\right)}
  }
}
\newcommand{\OOp}[2]{
  \ensuremath{
    \mathcal{O}
    \ifthenelse{\equal{#1}{0}}
      {}
      {\left({#1}\right)}
    \ifthenelse{\equal{#2}{1}}
      {}
      {^{\oplus{#2}}}
  }
}
\newcommand{\Proj}[1]{\ensuremath{\mathbb{P}^{#1}}}
\newcommand{\Sym}[2]{\ensuremath{\operatorname{Sym}^{#1}\left(#2\right)}}
\newcommand{\GL}[1]{\ensuremath{\operatorname{GL}\left(#1\right)}}
\newcommand{\gl}[1]{\mathfrak{gl}\left({#1}\right)}

\newcommand{\SL}[1]{\ensuremath{\operatorname{SL}\left(#1\right)}}
\newcommand{\PSL}[1]{\ensuremath{\mathbb{P}\operatorname{SL}\left(#1\right)}}
\newcommand{\Gr}[2]{\ensuremath{\operatorname{Gr}\left(#1,#2\right)}}
\newlength{\setBracketHeight}
\newcommand{\SetSuchThat}[2]{
  \settoheight{\setBracketHeight}{\ensuremath{#2}}
  \ensuremath{\left\{\left.{#1\rule{0cm}{\setBracketHeight}}\,
      \right|\,{#2}\right\}}}
\newlength{\tableRowHeight}
\newcommand{\tableExtraVerticalSpacing}{6pt}
\newcommand{\tablespacer}[1]
{
\settoheight{\tableRowHeight}{#1}
\addtolength{\tableRowHeight}{\tableExtraVerticalSpacing}
\ensuremath{\rule{0cm}{\tableRowHeight}#1}
}
\newcommand{\tablePMatrix}[1]
{
\tablespacer{\ensuremath{\begin{pmatrix}#1\end{pmatrix}}}
}
\newcommand{\tableBMatrix}[1]
{
\tablespacer{\ensuremath{\begin{bmatrix}#1\end{bmatrix}}}
}

\newcommand{\nForms}[2]{\ensuremath{\Omega^{#1} \left ( {#2} \right)}}
\newcommand{\Lm}[2]{\ensuremath{\Lambda^{#1}\left({#2}\right)}}
\newcommand{\Cohom}[2]{\ensuremath{H^{#1}\left({#2}\right)}}
\DeclareMathOperator{\Ad}{Ad}

\newcommand{\LieDer}{\ensuremath{\EuScript L}}
\newcommand{\hook}{\ensuremath{\mathbin{ \hbox{\vrule height1.4pt
        width4pt depth-1pt \vrule height4pt width0.4pt depth-1pt}}}}
\newcommand{\SO}[1]{\ensuremath{\operatorname{SO}\left(#1\right)}}
\newcommand{\Spin}[1]{\ensuremath{\operatorname{Spin}\left(#1\right)}}
\newcommand{\Symp}[1]{\ensuremath{\operatorname{Sp}\left(#1\right)}}
\newcounter{remarkCounter}
\newcommand{\pd}[2]{\frac{\partial #1}{\partial #2}}

\newcommand{\prol}[1]{\ensuremath{{#1}^{(1)}}}
\newcommand{\Gromov}[1]{\ensuremath{\check{#1}}}

\newcommand{\Bl}[2]{\ensuremath{\operatorname{Bl}_{#1}#2}}
\newcommand{\yes}{\checkmark}
\newcommand{\no}{{\textsf{x}}}

\newcommand{\ilim}{\mathop{\varprojlim}\limits}
\newcommand{\Thullen}{Thullen\xspace}
\newcommand{\cct}{complex codimension~$2$ or more\xspace}

\begin{document}

\allowdisplaybreaks

\renewcommand{\thefootnote}{$\star$}

\renewcommand{\PaperNumber}{058}

\FirstPageHeading

\ShortArticleName{Extension Phenomena for Holomorphic Geometric Structures}

\ArticleName{Extension Phenomena\\ for Holomorphic Geometric Structures\footnote{This paper is a
contribution to the Special Issue ``\'Elie Cartan and Dif\/ferential Geometry''. The
full collection is available at
\href{http://www.emis.de/journals/SIGMA/Cartan.html}{http://www.emis.de/journals/SIGMA/Cartan.html}}}

\Author{Benjamin MCKAY}

\AuthorNameForHeading{B.~McKay}

\Address{School of Mathematical Sciences, University College Cork, Cork, Ireland}
\Email{\href{mailto:b.mckay@ucc.ie}{b.mckay@ucc.ie}}
\URLaddress{\url{http://euclid.ucc.ie/pages/staff/Mckay}}

\ArticleDates{Received December 17, 2008, in f\/inal form May 07, 2009;  Published online June 08, 2009}

\Abstract{The most commonly encountered types of complex analytic $G$-structures and Cartan
geometries cannot have singularities of complex codimension 2 or more.}


\Classification{53A55; 53A20; 53A60; 53C10; 32D10}

{\tableofcontents}

\section{Introduction}

This article introduces methods to prove that various
holomorphic geometric structures (to be def\/ined below)
on complex manifolds cannot have singularities of \cct;
in other words, they extend holomorphically
across subsets of \cct.
Singularities can occur on complex hypersurfaces.
Simple (but perhaps less important) examples will show that
some other types of holomorphic
geometric structures \emph{can} have low dimensional
singularities.

Specif\/ically, we prove that
\begin{enumerate}\itemsep=0pt
 \item The underlying holomorphic principal bundle of a Cartan geometry
or $G$-structure extends holomorphically across a subset of
\cct
just when the Cartan geometry or $G$-structure does (see Theorem~\ref{theorem:CartanGeometry}
and Theorem~\vref{theorem:GeometricStructure}).
\item Holomorphic higher order structures extend across subsets of \cct
just when their underlying f\/irst order $G$-structures extend
(see Proposition~\vref{proposition:HigherOrder}).
\item The following holomorphic geometric structures extend holomorphically
across subsets of \cct:
\begin{enumerate}\itemsep=0pt
\item contact structures (see Theorem~\vref{theorem:ContactStructuresExtend}),
\item reductive $G$-structures (see Example~\vref{example:ReductiveAlgebraicG}),
\item scalar conservation laws (see Example~\vref{example:scalarConservationLaw}),
\item web geometries on surfaces (see Example~\vref{example:web}),
\item reductive Cartan geometries (see Example~\vref{example:reductiveAlgebraicCartan}) and
\item parabolic geometries (see Theorem~\vref{theorem:ParabolicGeometriesExtend}).
\end{enumerate}
\end{enumerate}
We also present a dictionary between holomorphic
extension problems for geometric structures
and holomorphic extension problems for maps
to complex homogeneous spaces.

The most striking example of our extension
theorems is that of 2nd order scalar ordinary
dif\/ferential equations (Example \vref{example:ODEs}),
for which we prove
that the geometry invariant under point
transformations extends across subsets of \cct,
and give examples in which the geometry invariant under
f\/iber-preserving transformations does not.

Every manifold and Lie group in these notes is
complex, and every map and bundle is holomorphic. I write
Lie groups as $G$, $H$, etc. and their Lie algebras
as $\mathfrak{g}$, $\mathfrak{h}$, etc.

\section{Def\/initions of Cartan geometries and f\/irst order structures}

For completeness, I will def\/ine the geometric structures of interest.

\subsection{First order structures}

\begin{Definition}
If $V \to M$ is a holomorphic vector bundle, the \emph{frame bundle}
of $V$, also called the \emph{associated principal bundle} and
denoted $FV$, is the set of complex linear
isomorphisms of f\/ibers of $V$ with some f\/ixed
vector space $V_0$. Clearly $FV \to M$ is a holomorphic
principal right $\GL{V_0}$-bundle, under the action
$r_h u=h^{-1}u$. When we need to be precise
about the choice of vector space $V_0$, we
will refer to $FV$ as the $V_0$-valued frame bundle.
\end{Definition}

\begin{Definition}
Suppose that $V_0$ is a complex vector space and $\rho : G \to \GL{V_0}$ is a
holomorphic representation of a complex Lie group $G$.
If $V \to M$ is a vector bundle and $FV$ its $V_0$-valued
frame bundle, then $G$ acts on $FV$ by having $g \in G$ take
$u \in FV \mapsto \rho(g)^{-1} u \in FV$.
We write this action as $r_g u = \rho(g)^{-1} u$.
A \emph{$G$-structure} on a manifold $M$ is a
principal $G$-bundle $E$ together with
a $G$-equivariant bundle morphism map $E \to FTM$.
See Gardner \cite{Gardner:1989} or Ivey and Landsberg
\cite{Ivey/Landsberg:2003} for an introduction to $G$-structures.

If we want to discuss $G$-structures for various
groups $G$, we will call them \emph{first order structures}.
\end{Definition}

\begin{Definition}
The \emph{standard flat $G$-structure} associated to a representation
$\rho : G \to \GL{V_0}$ is
the trivial bundle $E = V_0 \times G \to M = V_0$
mapped to $FV_0=V_0 \times \GL{V_0} \to V_0$
by $\left(v_0,g\right) \mapsto \left(v_0,\rho(g)\right)$.
\end{Definition}

\begin{Example}
Suppose that $E \to FTM$ is a $G$-structure for a complex representation
$G \to \GL{V_0}$. Let $K$ be the kernel of $G \to \GL{V_0}$,
and suppose that $G \to \GL{V_0}$ has closed image.
Then $E/K \subset FTM$ is a $G/K$-structure, called
the \emph{underlying embedded first order structure},
because it is embedded in $FTM$.
\end{Example}

\begin{Definition}
Fix a complex manifold $M$, a vector space $V_0$ with $\dim V_0=\dim M$,
and take $FTM$ the $V_0$-valued frame bundle.
Let $\pi : FTM \to M$ be the bundle mapping.
We def\/ine a~$V_0$-valued 1-form $\sigma$ on $FTM$, called the \emph{soldering form},
by $v \hook \sigma_{(m,u)} = u\left(\pi'(m,u)v\right)$,
for any $v \in T_{(m,u)} FTM$.
\end{Definition}

\subsection{Cartan geometries}

\begin{Definition}\label{def:CartanConnection}
Let $H \subset G$ be a closed subgroup of a Lie group,
with Lie algebras $\mathfrak{h} \subset \mathfrak{g}$.
A~\emph{Cartan geometry}
modelled on $G/H$ (also called a $G/H$-geometry)
on a manifold $M$ is a~choice of
principal right $H$-bundle $E \to M$, and 1-form $\omega \in
\nForms{1}{E} \otimes \mathfrak{g}$ called the \emph{Cartan connection},
which satisf\/ies the following conditions:
\begin{enumerate}\itemsep=0pt
\item
Denote the right action of $g \in G$ on $e \in E$ by $r_g e=eg$.
The Cartan connection transforms in the adjoint representation:
\[
r_g^* \omega = \Ad_g^{-1} \omega.
\]
\item
$\omega_e : T_e E \to \mathfrak{g}$ is a linear isomorphism at each point
$e \in E$.
\item
For each $A \in \mathfrak{g}$, def\/ine a vector f\/ield $\vec{A}$ on $E$ by
the equation $\vec{A} \hook \omega = A$. Then the vector f\/ields
$\vec{A}$ for $A \in \mathfrak{h}$ must generate the right $H$-action:
\[
\vec{A} =   \frac{d}{dt} r_{e^{tA}} \bigg|_{t=0}.
\]
\end{enumerate}
See Sharpe \cite{Sharpe:1997} for an introduction to Cartan
geometries.
\end{Definition}

\begin{Definition}
The \emph{standard flat $G/H$-geometry} is the Cartan
geometry whose bundle is $G \to G/H$ and whose Cartan
connection is $g^{-1} \, dg$.
\end{Definition}

\begin{Example}
Let $M$ be a complex manifold, and $\pi : E \to M$ a
$G/H$-geometry. Let $\mathfrak{g}$ and $\mathfrak{h}$
be the Lie algebras of $G$ and $H$. Let $V_0=\mathfrak{g}/\mathfrak{h}$.
Let $FTM$ be the $V_0$-valued frame bundle.
Let $\sigma=\omega+\mathfrak{h}$,
a semibasic 1-form valued in $V_0$.
At each point $e \in E$
the 1-form $\sigma$ determines a linear isomorphism
$u=u_e : T_m M \to V_0$
by the equation $u(\pi'(e)v) = v \hook \sigma$.
Map $e \in E \mapsto u_e \in FTM$. This map is an $H$-structure.
Let $H_1 \subset H$ be the subgroup of $H$
consisting of the elements of~$H$ which
act trivially on $V_0=\mathfrak{g}/\mathfrak{h}$.
The f\/ibers of the map $E \to FTM$
consist of the $H_1$-orbits in~$E$.
The map descends to a map $E/H_1 \to FTM$
called the \emph{underlying $\left(H/H_1\right)$-structure}
or \emph{underlying first order structure}
of the $G/H$-geometry.
\end{Example}

\begin{Definition}
The \emph{kernel} $K$ of a homogeneous space $G/H$ is
the largest normal subgroup of~$G$ contained in $H$:
\[
K = \bigcap_{g \in G} g H g^{-1}.
\]
\end{Definition}
As $G$-spaces, $G/H=(G/K)/(H/K)$, so strictly speaking, the kernel is not an invariant
of the homogeneous space, but of the choice of
Lie group $G$ and closed subgroup $H$, but
we will ignore this subtlety.

\begin{Lemma}[Sharpe \cite{Sharpe:1997}]
Suppose that $K$ is the kernel of $G/H$, with Lie
algebra $\mathfrak{k}$, and that $E \to M$ is a Cartan geometry modelled
on $G/H$, with Cartan connection $\omega$. Let
$E'=E/K$ and $\omega'=\omega + \mathfrak{k} \in \nForms{1}{E}
\otimes \left(\mathfrak{g}/\mathfrak{k}\right)$. Then $\omega'$ drops to a
$1$-form on $E'$, and $E' \to M$ is a Cartan geometry, called the
\emph{reduction} of $E$.
\end{Lemma}

\begin{Definition}
A complex homogeneous space $G/H$ is
called \emph{reduced} or \emph{faithful} if $K=1$,
and \emph{almost reduced} or \emph{almost faithful}
if $K$ is a discrete subgroup of $H$.
\end{Definition}

\subsection{Parabolic geometries}

\begin{Definition}
A Lie subalgebra $\mathfrak{p} \subset \mathfrak{g}$ of a semisimple Lie algebra
is called \emph{parabolic} if it contains a Borel subalgebra (see Serre \cite{Serre:2001}). A
connected Lie subgroup  of a semisimple Lie group is called
\emph{parabolic} if its Lie algebra is parabolic.
A homogeneous space $G/P$
(with $P$ parabolic and $G$ connected) is called a
\emph{rational homogeneous variety}; see
Landsberg \cite{Landsberg:2005} for an introduction
to rational homogeneous varieties.
\end{Definition}
\begin{Definition}\label{def:parabolicGeometry}
A \emph{parabolic geometry} is a Cartan geometry modelled on a
rational homogeneous variety. See \v{C}ap \cite{Cap:2005}
for an introduction to parabolic geometries.
\end{Definition}

\subsection{Curvature}

\begin{Definition}
The \emph{curvature} of a Cartan geometry
is the 2-form $d \omega + \frac{1}{2} \left[\omega,\omega\right]$.
A Cartan geometry is \emph{flat} if its curvature vanishes.
\end{Definition}

\subsection{Development}

\begin{Definition}
Suppose that $E_0 \to M_0$ and $E_1 \to M_1$ are two
$G/H$-geometries, with Cartan connections $\omega_0$ and $\omega_1$,
and $X$ is manifold, perhaps with boundary and corners. A smooth map
$\phi_1 : X \to M_1$ is a \emph{development} of a smooth map $\phi_0
: X \to M_0$ if there exists a smooth isomorphism $\Phi : \phi_0^*
E_0 \to \phi_1^* E_1$ of principal $H$-bundles
so that $\Phi^* \omega_1 = \omega_0$.
\end{Definition}

\begin{Definition}
Suppose that $E_0 \to M_0$ and $E_1 \to M_1$ are $G/H$-geometries.
Suppose that $\phi_1 : X \to M_1$ is a development of a smooth map
$\phi_0 : X \to M_0$ with isomorphism $\Phi : \phi_0^* E_0 \to
\phi_1^* E_1$. By analogy with Cartan's method of the moving frame,
we will call $e_0$ and $e_1$ \emph{corresponding frames} of the
development if $\Phi\left(e_0\right)=e_1$.
\end{Definition}

\begin{Lemma}\label{lemma:developCurves}
Suppose that $C$ is a simply connected Riemann surface,
that $E_0 \to M_0$ is a~$G/H$-geometry,
that $\phi_0 : C \to M_0$, and that $e_0 \in \phi_0^* E_0$.
Then $\phi_0$ has a unique development $\phi_1 : C \to G/H$
to the model with a unique isomorphism $\Phi : \phi_0^* E_0 \to \phi_1^* G$
so that $\Phi\left(e_0\right)=1 \in G$.
\end{Lemma}
\begin{proof} (This proof is adapted from
McKay \cite{McKay:2007}.)
The local existence and
uniqueness of a~development is clear by applying the Frobenius
theorem to the Pfaf\/f\/ian system $g^{-1} \, dg = \omega_0$ on $\phi_0^*
E_0 \times G$. (The curvature does not af\/fect the involutivity
of this Pfaf\/f\/ian system.) The maxi\-mal connected integral manifolds project
locally dif\/feomorphically to $\phi_0^* E_0$, because $\omega_0$ is a
coframing on them.

Because $C$ is simply connected,
$\phi_0^*  E_0 \to C$ is a trivial bundle, with a global section
$s_0$.  If we can develop, then this global section is identif\/ied
via the isomorphism $\Phi$ with a global section $s_1 : C \to
\phi_1^* G$ so that
\begin{equation}\label{eqn:devel}
s_1^* g^{-1} \, dg = s_0^* \omega_0.
\end{equation}
Conversely, if we can solve this equation, then there is a unique
isomorphism $\Phi$ for which
\[
\Phi\left(s_0 h\right) = s_1 h
\]
for all $h \in H$, by triviality of the bundles. So it suf\/f\/ices to
solve equation~\eqref{eqn:devel}.

Equation~\eqref{eqn:devel} is an ordinary
dif\/ferential equation of Lie type,
\[
g^{-1}  dg = s_0^* \omega_0,
\]
(see Bryant \cite[p.~55]{Bryant:1995}
or Sharpe \cite[p.~118]{Sharpe:2002})
so has a unique global solution with given initial condition $g=g_0$
at $t=t_0$. For example, if $G$ is a Lie subgroup of $\GL{n,\R{}}$
for some $n$, then global existence and uniqueness of a development
follow from writing the ordinary dif\/ferential equations of Lie type
as a linear ordinary dif\/ferential equation:
\[
dg = g  s_0^* \omega_0.
\]
More generally, one glues together local solutions by using the
group action of $G$ to make two local solutions match up at some
point. Then the local solutions match near this point by uniqueness.
Any compact simply connected subset of $C$ is covered by f\/initely
many domains of such local solutions, which thereby must patch to a
global solution.
\end{proof}

\begin{Lemma}
 Suppose that $E_0 \to M_0$ is a flat $G/H$-geometry,
that $X$ is a simply connected complex manifold,
that $\phi_0 : X \to M_0$ is a holomorphic map, and
that $e_0 \in \phi_0^* E_0$. Then there is a unique
developing map $\phi_1 : X \to G/H$ with isomorphism
$\Phi : \phi_0^* E_0 \to \phi_1^* G$ for which
$\Phi\left(e_0\right)=1 \in G$.
\end{Lemma}
\begin{proof}
 Locally, this is the Frobenius theorem. Just as for
curves above, the local developments patch together
under $G$-action to extend uniquely to all of $X$.
\end{proof}

\section{Extensions of maps}

\subsection{Def\/initions and examples of extension phenomena}

\begin{Definition}
 A subset $S \subset M$ of a complex manifold
is of \emph{\cct} if each point $s \in S$
lies in an open set $U_s \subset M$
so that $S \cap U_s$ is contained in
a complex analytic subvariety $V_s$ of $U_s$, of
complex codimension 2 or more.
\end{Definition}

\begin{Lemma}[Hartogs extension lemma~\cite{Hartogs:1906}]
If $M$ is a complex manifold $M$ and $S \subset M$
is of \cct, then every holomorphic function
$f : M \backslash S \to \C{}$ extends
to a unique holomorphic function
$f : M \to \C{}$.
\end{Lemma}
We will paraphrase the Hartogs extension lemma as
saying that holomorphic functions
extend over subsets of \cct.
Krantz \cite{Krantz:2008} provides an introduction to Hartogs extension phenomena,
while Merker and Porten \cite{Merker/Porten:2007} give an elegant
new proof of the result in a~more general form.

Consider two extension problems for holomorphic
functions, holomorphic bundles, or holomorphic
geometric structures:
\begin{itemize}\itemsep=0pt
 \item
the \emph{Hartogs extension problem} of extending
holomorphically from a domain in a Stein manifold
to its envelope of holomorphy
\item
the \emph{\Thullen extension problem} of extending
holomorphically across a subset of \cct.
\end{itemize}
\begin{Remark}
The expression \emph{Thullen-type extension}
is used in complex analysis to mean extension
of a holomorphic function
(or extension of a holomorphic section of a vector bundle,
or extension of a holomorphic vector bundle, etc.)
from $M \backslash S$ to $M$, where $M$
is a complex manifold, and $S \subset M$
has the form $S = H \backslash U$,
where $H \subset M$ is a complex hypersurface,
and $U \subset M$ is an open set intersecting
every codimension 1 component of $H$.
In this paper, we will use the term
\emph{\Thullen extension} only in the narrower
sense above, for simplicity.
\end{Remark}

\begin{Definition}
A complex space $X$ is a \emph{Hartogs extension target}
if every holomorphic
map $f : D \to X$ from a domain $D$ in a Stein
manifold extends to the envelope of holomorphy of $D$.
\end{Definition}
We will also need the following two
simpler and weaker properties, for \Thullen
extension problems.
\begin{Definition}
We will say that a complex space $X$ is a \emph{\Thullen extension target}, to mean
that every holomorphic map $f : M \backslash S \to X$
extends to a holomorphic map $f : M \to X$, where
$M$ is any complex manifold and $S \subset M$
is any subset of \cct.
(Informally, we will also say that
holomorphic maps to $X$ extend across
subsets of \cct.)
Hartogs extension targets are \Thullen extension targets.
Similarly, we will say that $X$ is a \Thullen extension
target for local biholomorphisms to mean that local
biholomorphisms to $X$ extend across
subsets of \cct, etc.
\end{Definition}

\begin{Example}
$\C{}$ is a Hartogs extension target.
\end{Example}

\begin{Example}
If $X$ is any complex manifold of dimension at least two,
and $x \in X$, then $X \backslash x$
is neither a Hartogs extension target nor a \Thullen extension target for
local biholomorphisms: $\text{id} : X \setminus x \to X \setminus x$
doesn't extend to a map $X \to X \setminus x$.
\end{Example}

\begin{Example}
The map $z \in \C{n+1} \backslash 0 \to \C{}z \in \Proj{n}$
doesn't extend over the puncture at $0$.
\end{Example}

\begin{Example}
Let $\Bl{m}{M}$ be the blowup of a complex
manifold $M$ at a point $m$.
Map $M \backslash m \to \Bl{m}{M}$ by the obvious local biholomorphism.
This local biholomorphism clearly doesn't extend across the puncture.
\end{Example}

\begin{Example}\label{example:extensionSubspaces}
Every closed complex subspace of a~Hartogs/\Thullen extension target
is a~Hartogs/\Thullen extension target.
\end{Example}

\begin{Example}
 If $X$ is the blowup of a complex manifold along
an analytic subvariety, then $X$ is \emph{not}
a Hartogs or \Thullen extension target (Hartogs and
\Thullen extension targets are \emph{minimal}).
\end{Example}

\begin{Example}[Ivashkovich \protect{\cite[Proposition 3, p.~195]{Ivashkovich:1985}}]
A product is a Hartogs/\Thullen extension target just when the factors are.
\end{Example}

\begin{Example}
Af\/f\/ine analytic varieties and Stein manifolds are
Hartogs and \Thullen extension targets. (Recall
that every Stein manifold admits a closed holomorphic
embedding into $\C{N}$ for some integer $N>0$.)
\end{Example}

\begin{Example}
Pseudoconvex domains in Hartogs/\Thullen extension
targets are themselves Hartogs/\Thullen extension targets.
\end{Example}

\begin{Example}
Let $S$ be the Hopf surface: $\left(\C{2} \backslash 0\right)/(z \sim 2z)$.
The Hopf surface is a compact complex surface.
Take the map $f : \C{2} \backslash 0 \to S$ taking
each point $z$ to its equivalence class $[z] \in S$.
This map is a local biholomorphism onto $S$, but
doesn't extend to $\C{2}$ because distinct points of~$S$
have preimages arbitrarily close to $0$.
Therefore the Hopf surface is not an extension
target in either sense.
See Wehler~\cite{Wehler:1982} for an introduction to Hopf surfaces.
\end{Example}

\begin{Lemma}[Ivashkovich \protect{\cite[Lemma 6, p.~229]{Ivashkovich:1986}}]\label{lemma:coverings}
If $X$ and $Y$ are complex manifolds
and $X \to Y$ is an unramified covering map,
then $X$ is a Hartogs/\Thullen extension target
just when~$Y$~is.
\end{Lemma}
\begin{proof}
Consider the \Thullen extension problem
of extending a map $M \setminus S \to X$.
Local extensions will clearly glue together
to give a global extension. So we can
replace our manifold $M$ with a simply
connected open subset of $M$. Take
$S \subset M$ any subset of \cct.
Then $M \backslash S$ is also simply connected, so we
can replace $X$ by any covering space of~$X$.
The Hartogs extension problem is more subtle;
see Ivashkovich~\cite{Ivashkovich:1986}.
\end{proof}

\begin{Lemma}\label{lemma:fundamentalGroupEnvelope}
If $M$ is a domain in Stein
manifold, and $\hat{M}$ is the
envelope of holomorphy of~$M$,
then $\pi_1\left(M\right) \to \pi_1\left(\hat{M}\right)$
is surjective.
\end{Lemma}
\begin{proof}
Let $\Gamma$ be the image of
the morphism $\pi_1\left(M\right) \to \pi_1\left(\hat{M}\right)$
given by the inclusion $M \subset \hat{M}$.
Let $\tilde{M}$ be the universal
covering space of $\hat{M}$,
and let $M'=\tilde{M}/\Gamma$,
with projection mapping $\pi : M' \to \hat{M}$.
Clearly there is a section $M \to M'$
over $M$. By Lemma~\vref{lemma:coverings},
since $M'$ is a~covering space of $\hat{M}$,
the holomorphic map $M \to M'$
extends uniquely to a holomorphic
map $\hat{M} \to M'$. By
analytic continuation, this
map must be a section of $M' \to \hat{M}$,
so $M'=\hat{M}$.
\end{proof}

\begin{Lemma}
If a local biholomorphism extends across a subset of \cct to
a holomorphic map, then it extends uniquely to a local biholomorphism.
\end{Lemma}
\begin{proof}
Express the extended map in some local coordinates
as a holomorphic function $w(z)$, $z$,~$w$~points of some open subsets of $\C{n}$.
Then $\det w'(z) \ne 0$ away from $z \in S$. But then $\det w'(z)$ and
$1/\det w'(z)$ extend holomorphically across $S$ by Hartogs extension lemma.
So $\det w'(z) \ne 0$ throughout the domain of the local coordinates.
\end{proof}

\begin{Lemma}
If $X$ is a Hartogs/\Thullen extension target $($for local biholomorphisms$)$,
and $F_{\alpha} : X \to \C{}$
are some holomorphic functions,
then $X \backslash \cup_{\alpha} \left(F_{\alpha}=0\right)$
is a Hartogs/\Thullen extension target $($for local biholomorphisms$)$.
\end{Lemma}
\begin{proof}
Suppose that
$f : M \backslash S \to X \backslash \cup_{\alpha} \left(F_{\alpha}=0\right)$.
Then $f$ extends to a map $f : M \to X$. Away
from $S$, $1/\left(F_{\alpha} f\right)$ is a holomorphic function, so
extends to $M$, and therefore $F_{\alpha} \ne 0$ on $M$.
\end{proof}

\begin{Example}
If $f : M \backslash S \to G$ and $G \subset \GL{n,\C{}}$
is a closed Lie subgroup, then $f$ and $f^{-1}$  extend
holomorphically to matrix-valued functions, so
$f$ extends to a holomorphic map to $\GL{n,\C{}}$
and therefore $f$ extends to a map to $G$. So closed complex subgroups
of $\GL{n,\C{}}$ are \Thullen extension targets. Similarly,
they are Hartogs extension targets. So are their
covering spaces, by Lemma~\vref{lemma:coverings},
so any Lie group which admits a faithful complex
Lie algebra representation is a Hartogs
and \Thullen extension target.
\end{Example}

\begin{Lemma}[Adachi, Suzuki, Yoshida \cite{Adachi/Suzuki/Yoshida:1973}]\label{lemma:ComplexLieGroups}
Complex Lie groups are Hartogs and \Thullen extension targets.
\end{Lemma}
We will f\/irst give a short proof
that they are \Thullen extension targets.
\begin{proof}
Suppose that $M$ is a complex manifold,
$S \subset M$ is a subset of \cct and
$f : M \backslash S \to G$
is a holomorphic map to a complex Lie group.
We can assume that $M$ and $G$ are connected.
Local extensions will glue together to
produce a global extension, so we can
assume that $M$ is simply connected, and therefore
that $M \backslash S$ is too.
The 1-form $f^* g^{-1} \, dg$ extends to
a 1-form on $M$: just write it out
in local coordinates near $S$ and use Hartogs
extension lemma. Now apply the fundamental
theorem of calculus for maps to Lie groups:
Sharpe \cite[p.~124]{Sharpe:1997}.
\end{proof}

To carry solve Hartogs extension problems, we will
need to extend meromorphic functions as well.
\begin{Theorem}[Kajiwara and Sakai \protect{\cite[p.~75]{Kajiwara/Sakai:1967}}]\label{theorem:KajiwaraSakai}
All meromorphic functions on a domain~$M$ in a Stein manifold
extend meromorphically to the envelope of holomorphy of~$M$.
\end{Theorem}

\begin{Lemma}
Any holomorphic vector field on a domain in
a Stein manifold extends holomorphically
to the envelope of holomorphy.
\end{Lemma}
\begin{proof}
Suppose that $M$ is a domain in a Stein
manifold, and that $\hat{M}$ is the
envelope of holomorphy of $M$; in
particular $\hat{M}$ is Stein,
and all holomorphic functions on $M$
extend to $\hat{M}$ (see H\"ormander \cite[p.~116]{Hormander:1990}).

To extend a vector f\/ield $v$ from $M$ to
$\hat{M}$, take a holomorphic
function $h : \hat{M} \to \C{}$, and extend the holomorphic
function $\LieDer_v h$ to $\hat{M}$. The Leibnitz identity
extends. For each point of a Stein manifold $\hat{M}$,
there is a holomorphic map $\hat{M} \to \C{n}$
which is a local biholomorphism near that point
(by def\/inition; see
H\"ormander \cite[p.~116]{Hormander:1990}).
Therefore a vector f\/ield is uniquely
determined by its action on holomorphic
functions.
\end{proof}

\begin{Lemma}\label{lemma:hypersurfaceIntersections}
Suppose that $M$ is a domain in a Stein
manifold. Let $\hat{M}$ be the envelope
of holomorphy of $M$. Then $M$ intersects every
complex hypersurface in $\hat{M}$.
\end{Lemma}
\begin{proof}
Let $H \subset \hat{M}$
be a complex hypersurface.
Henri Cartan's Theorem A
(H\"ormander \cite[Theorem~7.2.8, p.~190]{Hormander:1990})
says that for every coherent
sheaf $F$ on any Stein manifold
every f\/iber~$F_m$ is generated by
global sections.
Let $\OO{H}$ be the
line bundle corresponding
to the hypersurface~$H$.
By Cartan's theorem, for each point $m \in H$,
there must be a global section
of $\OO{H}$ not vanishing at $m$.
Take such a section, and
let $H'$ be its zero locus.

Henri Cartan's theorem B
(H\"ormander \cite[Theorem~7.4.3, p.~199]{Hormander:1990})
says that the cohomology of
any coherent sheaf
on any Stein manifold vanishes
in positive degrees.
By the exponential sheaf sequence,
(see H\"ormander \cite[p.~201]{Hormander:1990})
$\Cohom{1}{\hat{M},{\mathcal{O}}^{\times}}=\Cohom{2}{\hat{M},\Z{}}$.
In other words, line bundles on $\hat{M}$ are precisely
determined by their coholomogy
classes. Therefore a divisor
is the divisor of a meromorphic
function just when it has trivial
f\/irst Chern class. In particular,
$H'-H$ is the divisor of a meromorphic
function, say $f : \hat{M} \to \C{}$.
So $f$ is holomorphic on $\hat{M} \setminus H$,
and has poles precisely on $H$.
If $H$ does not intersect $M$, then
$f$ is holomorphic on $M$,
so extends holomorphically
to $\hat{M}$, i.e.\ has no poles, so $H$ is empty.
\end{proof}

\begin{Corollary}\label{corollary:meromorphicHolomorphic}
Suppose that $M$ is a domain in a Stein
manifold, with envelope of holomorphy
$\hat{M}$ and that $H \subset \hat{M}$
is a complex hypersurface. Every
function meromorphic on $M$ and
holomorphic on $M \setminus H$
extends uniquely to a function
meromorphic on $\hat{M}$ and
holomorphic on~$\hat{M} \setminus H$.
\end{Corollary}
\begin{proof}
Suppose that $f$ is meromorphic
on $M$ and holomorphic on $M \setminus H$,
so the poles of $f$ lie on $H$.
By the theorem of
Kajiwara and Sakai (Theorem~\vref{theorem:KajiwaraSakai}),
every meromorphic function on a
domain $M$ in a Stein manifold extends
to a meromorphic function on $\hat{M}$.
Suppose that $f$ extends to have
poles on some hypersurface $H \cup Z$.
We can assume that $Z$ contains no
component of $H$, so that $H \cap Z$ is
a subset of \cct.
But then every component of $Z$ must be
a hypersurface not intersecting $M$,
so must be empty.
\end{proof}

\begin{Lemma}[H\"ormander \cite{Hormander:1990} p. 116]
Suppose that $M$ is a domain in an
$n$-dimensional Stein
manifold, with envelope of holomorphy
$\hat{M}$.
For each point $m \in \hat{M}$,
there is a map $f : \hat{M} \to \C{n}$
which is a local biholomorphism near $m$.
\end{Lemma}

\begin{Proposition}
Let $M$ be a domain in Stein
manifold, and $\hat{M}$ be
the envelope of holomorphy
of $M$.
Every holomorphic tensor extends
from $M$ to $\hat{M}$.
\end{Proposition}
\begin{Remark}
This proposition generalizes a theorem of Aeppli \cite{Aeppli:1965}.
\end{Remark}
\begin{proof}
For simplicity, suppose that
we wish to extend a holomorphic 1-form
$\omega$ on $M$ to $\hat{M}$.
Pick a point $m_0 \in \hat{M}$.
Pick a map $f : \hat{M} \to \C{n}$
which is a local biholomorphism
near $m_0$, say
on some open set $\hat{M} \setminus H_f$
containing $m_0$.
Write the 1-form $\omega$ on $M$ as $\omega=g \, df$.
Then $g$ is a~meromorphic function
on $M$, and holomorphic
on $M \setminus H_f$.
By Corollary~\vref{corollary:meromorphicHolomorphic},
$g$ extends uniquely to a~function
meromorphic
on $\hat{M}$ and holomorphic
on $\hat{M} \setminus H_f$.
So we can extend $\omega$
to a neighborhood of $m_0$
by $\omega = g \, df$.
Clearly the extension is
unique where def\/ined.
Since $m_0$ is an arbitrary
point of $\hat{M}$, $\omega$
extends holomorphically to all
of $\hat{M}$.
\end{proof}

We now f\/inish the proof of Lemma~\vref{lemma:ComplexLieGroups}:
complex Lie groups are Hartogs and Thullen extension targets.
\begin{proof}
Suppose that $M$ is a domain in a Stein
manifold, with envelope of holomorphy $\hat{M}$,
and $f : M \to G$
is a holomorphic map to a complex Lie group.
We can assume that $M$ and $G$ are connected.
The 1-form $f^* g^{-1} \, dg$ extends to
a 1-form on $\hat{M}$. Apply the fundamental
theorem of calculus for maps to Lie groups:
Sharpe \cite[p.~124]{Sharpe:1997}
to see that some covering space
$\pi : \tilde{M}$ of $\hat{M}$ bears a holomorphic
map $\tilde{f} : \tilde{M} \to G$ satisfying
$\tilde{f}^* g^{-1} \, dg=\omega$ (where $\omega$
here denotes the pullback of $\omega$ on
$\hat{M}$ to the covering space $\tilde{M}$).
Note that we don't require $\tilde{M}$ to
be the universal covering space.
Moreover,
there is a morphism of groups $h : \pi_1\left(\hat{M}\right) \to G$
so that $\tilde{f}\left(\gamma m\right) = h(\gamma) \, \tilde{f}(m)$
for any $\gamma \in \pi_{1}\left(\hat{M}\right)$
and $m \in \tilde{M}$. This map $\tilde{f}$ is uniquely
determined up to the action of $G$.
Pick a point $\tilde{m}_0 \in \tilde{M}$ so that the
point $m_0 = \pi\left(\tilde{m}_0\right)$ lies in $M$.
Arrange by $G$-action
that $\tilde{f}\left(\tilde{m}_0\right)=f\left(m_0\right)$.
Then clearly
$\tilde{f}\left(\tilde{m}\right)=f\left(m\right)$
whenever $m=\pi\left(\tilde{m}\right) \in M$.
In particular, $h=1$ on $\pi_1(M)$.
By Lemma~\ref{lemma:fundamentalGroupEnvelope}, $h=1$ on $\pi_1\left(\tilde{M}\right)$.
\end{proof}

The three following results of Matsushima and Morimoto provide
a larger class of examples of Stein manifolds.

\begin{Lemma}[Matsushima and Morimoto \protect{\cite[Theorem~2, p.~146]{Matsushima/Morimoto:1960}}]\label{lemma:MatMoriGroup}
A complex Lie group $G$ is a Stein manifold
just when the  identity component of
the center of $G$ is a linear
algebraic group, i.e.\ contains no complex tori.
\end{Lemma}

\begin{Lemma}[Matsushima and Morimoto \protect{\cite[Theorem~4, p.~147]{Matsushima/Morimoto:1960}}]\label{lemma:MatMori}
If a holomorphic principal bundle
has Stein base, and Stein fibers,
then it has Stein total space.
\end{Lemma}

\begin{Lemma}[Matsushima and Morimoto \protect{\cite[Theorem~6, p.~153]{Matsushima/Morimoto:1960}}]\label{lemma:MatMoriTwo}
If a holomorphic principal bundle $G$-bundle $E \to M$,
has Stein base $M$, and if $H \subset G$ is a
closed connected complex Lie subgroup,
and $G/H$ is Stein,
then $E/H \to M$ has Stein total space.
\end{Lemma}

\begin{Lemma}\label{lemma:ClassicalInvariantTheory}
Complex reductive homogeneous spaces $($i.e. $G/H$ with $H \subset G$
a closed reductive algebraic subgroup of a complex
linear algebraic group $G)$ are affine analytic varieties.
\end{Lemma}
Therefore reductive homogeneous spaces are Hartogs and \Thullen extension targets.
\begin{proof}
The categorical quotient $G//H$ (spectrum of $H$-invariant polynomials)
of af\/f\/ine algebraic varieties by reductive algebraic groups
are af\/f\/ine algebraic varieties;
see Procesi \cite[Theo\-rem~2, p.~556]{Procesi:2007}.
The quotient will in general only parameterize the closed
orbits, and admit a holomorphic submersion $G/H \to G//H$.
But the $H$-orbits on $G$ are
the translates of $H$, so all orbits are closed, and thus $G/H=G//H$.
\end{proof}

\begin{Corollary}\label{corollary:ExtendReductiveReduction}
Suppose that $G$ is a complex linear algebraic group
and $H \subset G$ is a closed complex
reductive algebraic Lie subgroup.
Suppose that $M$ is a domain in a Stein
manifold, and that $\hat{M}$ is the
envelope of holomorphy of $M$.
If $E \to \hat{M}$ is a holomorphic
principal $G$-bundle,
then every holomorphic section of $\left.E/H\right|_{M} \to M$
extends uniquely to a holomorphic
section of $E/H \to \hat{M}$.
\end{Corollary}
\begin{proof}
Let $H^0 \subset H$ be the connected
component of the identity element.
Then $E/H^0$ is a Stein manifold,
by Lemma~\vref{lemma:MatMoriTwo},
so a Hartogs extension target.
By Lemma~\vref{lemma:coverings}, $E/H$
is a Hartogs extension target.
Therefore any section
$s$ of $\left.E/H\right|_{M} \to M$
extends to a unique map $s : \hat{M} \to E/H$.
Consider the bundle map $\pi : E \to \hat{M}$.
Over $M$, $s$ satisf\/ies $\pi s = \text{id}$.
By analytic continuation, this holds over
$\hat{M}$ as well.
\end{proof}

\begin{Lemma}
Suppose that $X$ is a \Thullen extension target for local biholomorphisms.
The complement of any hypersurface in $X$ is also a \Thullen
extension target for local biholomorphisms.
\end{Lemma}
\begin{proof}
If $f : M \backslash S \to X \backslash H$ is a local biholomorphism,
then it extends to a local biholomorphism $f : M \to X$. The hypersurface
$f^{-1}(H)$ is either empty or else cannot be contained in $S$, and
so intersects $M \backslash S$, so $f$ doesn't take $M \backslash S$
to $X \backslash H$.
\end{proof}

\begin{Lemma}
Suppose that $X$ is a Hartogs extension target for local biholomorphisms.
The complement of any complex hypersurface in $X$ is also a Hartogs
extension target for local biholomorphisms.
\end{Lemma}
\begin{proof}
Suppose that $M$ is a domain in a Stein manifold, with
envelope of holomorphy $\hat{M}$.
If $f : M  \to X \backslash H$ is a local biholomorphism,
then it extends to a local biholomorphism $f : \hat{M} \to X$.
The complex hypersurface $f^{-1}(H)$ cannot intersect $M$,
so is empty by Lemma~\vref{lemma:hypersurfaceIntersections}.
\end{proof}

\begin{Proposition}\label{proposition:VectorBundleSections}
Suppose that $M$ is a domain in a Stein
manifold with envelope of holomorphy~$\hat{M}$.
If $V \to \hat{M}$ is a holomorphic vector bundle,
then every section of $V$ over $M$
extends to a section of $V$ over $\hat{M}$.
\end{Proposition}
\begin{proof}
Pick a section $f$ of $V$ over $M$,
and a point $\hat{m} \in \hat{M}$.
By Cartan's theorem A, there are
global sections $f_1, f_2, \dots, f_p$ of
$V$ over $\hat{M}$ whose values
at $\hat{m}$ are a basis of $V_{\hat{m}}$.
These sections are linearly independent
away from some hypersurface $H \subset \hat{M}$.
On $\hat{M} \setminus H$, $V$ is trivial.
So on~$M \setminus H$,
we can write $f = \sum_i c_i f_i$ for
some holomorphic functions
$c_i$ on $M \setminus H$.

Check that the functions $c_i$
are meromorphic on $M$ as follows.
In any local trivialization $g_1, g_2, \dots, g_p$
for $V$, write
$f_i=\sum_j A_{ij} g_j$ and $f=\sum_j b_j g_j$.
Let $A = \left(A_{ij}\right)$.
The entries of  $A^{-1}$ are rational functions
of the entries of $A$, so we can
write $f=\sum_{ij} A^{-1}_{ij} b_i f_i$.

Each $c_i$ extends to a holomorphic
function on $\hat{M} \setminus H$,
and therefore we can extend $f$ to
$f= \sum_i c_i f_i$. This extends
$f$ to a neighborhood of the arbitrary
point $\hat{m}$.
\end{proof}

\begin{Corollary}\label{lemma:holomorphicConnectionsExtend}
Suppose that $M$ is a domain in a Stein
manifold with envelope of holomorphy~$\hat{M}$. If $E \to \hat{M}$
is a holomorphic $G$-bundle,
then every holomorphic connection
on $\left.E\right|_M$ extends
uniquely to a holomorphic
connection on $E$.
\end{Corollary}
\begin{proof}
Suppose that $\pi : E \to \hat{M}$
is a principal $G$-bundle,
where $G$ is a complex Lie
group with Lie algebra $\mathfrak{g}$.
Let $E_0 = \left.E\right|_{M}$.
A connection on $E_0 \to M$ is a section
of the vector bundle
$\left(T^*E_0 \otimes \mathfrak{g}\right)^G \to M$,
i.e.\ the vector bundle whose local sections are
precisely the $G$-invariant
local sections of $T^*E_0 \otimes \mathfrak{g} \to E_0$. This
vector bundle extends to the
holomorphic vector bundle
$\left(T^*E \otimes \mathfrak{g}\right)^G \to \hat{M}$.
By Proposition~\ref{proposition:VectorBundleSections},
every section of $\left(T^*E_0 \otimes \mathfrak{g}\right)^G \to M$
extends to a~holomorphic
section of
$\left(T^*E \otimes \mathfrak{g}\right)^G \to \hat{M}$.
So we can extend $\omega$ to a $\mathfrak{g}$-valued
1-form on $E$.

For each vector $A \in \mathfrak{g}$,
denote by $\vec{A}$ the associated
vertical vector f\/ield on $E$ giving
the inf\/initesimal $\mathfrak{g}$-action.
A connection on $E \to \hat{M}$
is precisely a $G$-equivariant section $\omega$ of
$T^*E_0 \otimes \mathfrak{g} \to M$,
for which $\vec{A} \hook \omega=A$,
for $A \in \mathfrak{g}$.
Clearly this identity must extend
from $E_0$ to $E$ by analytic continuation.
\end{proof}

\begin{Lemma}\label{lemma:SteinHasHolomorphicConnections}
Any holomorphic principal bundle on any
Stein manifold admits a holomorphic connection.
\end{Lemma}
\begin{proof}
Let $G$ denote a complex Lie group.
The obstruction to a holomorphic
connection on a~holomorphic principal
$G$-bundle $\pi : E \to M$
is the Atiyah class
\[
a\left(M,E\right) \in \Cohom{1}{M,T^*M \otimes \Ad(E)};
\]
see Atiyah \cite{Atiyah:1957}.
Henri Cartan's theorem B says that all
positive degree cohomology groups
of coherent sheaves on Stein manifolds
vanish;
see H\"ormander \cite[Theorem~7.4.3, p.~199]{Hormander:1990}.
Therefore the Atiyah class
vanishes, and so there is a holomorphic connection.
\end{proof}

\begin{Corollary}
If $M$ is a domain in a Stein manifold,
and $E \to M$ is a holomorphic principal bundle
with nonzero Atiyah class, then
$E$ does not extend to a holomorphic
principal bundle on the envelope of
holomorphy of $M$.
\end{Corollary}

\begin{Conjecture}
A holomorphic principal bundle $E \to M$ over
a domain $M$ in a Stein manifold
extends uniquely to a holomorphic
principal bundle over the envelope
of holomorphy~$\hat{M}$ just when
the Atiyah class of $E$ vanishes.
\end{Conjecture}

\subsection{Meromorphic extension theorems}

\begin{Definition}
 If $X$ and $Y$ are complex manifolds, and $X$ is connected,
a \emph{meromorphic map}
(in the sense of Remmert \cite[Def\/inition~15, p.~367]{Remmert:1957}) $f : X \to Y$ is a
choice of nonempty compact set $f(x) \subset Y$ for each $x \in X$,
so that
\begin{enumerate}\itemsep=0pt
 \item $f(x)$ is a single point for $x$ in some dense open subset of $X$,
 \item the pairs of $(x,y)$ with $y \in f(x)$ form
an irreducible analytic variety $\Gamma \subset X \times Y$,
called the \emph{graph} of $f$.
\end{enumerate}
\end{Definition}

Clearly the composition $fg$ of a holomorphic map $f$ with a meromorphic
map $g$, given by mapping points $(x,y)$ of the graph to $(x,f(y))$, is
also a meromorphic map. It is easy to prove
(see Remmert~\cite{Remmert:1957},
Siu~\cite{Siu:1974}, Ivashkovich~\cite{Ivashkovich:1992}) that there is
a subset of \cct in $X$ (called the
\emph{indeterminacy locus}) over
which $f$ is the graph of a holomorphic map.

The \Thullen extension problem:
\begin{Theorem}[Siu \cite{Siu:1974}]\label{theorem:Siu}
Take a complex manifold $M$ and a subset of \cct $S \subset M$.
Every holomorphic map $f : M \backslash S \to X$
to a compact K\"ahler manifold $X$ extends to a meromorphic map $f : M \to X$.
\end{Theorem}

The Hartogs extension problem:
\begin{Theorem}[Ivashkovich \cite{Ivashkovich:1992}]\label{theorem:Ivashkovich}
 Every holomorphic map from a domain in a Stein manifold
to a compact K\"ahler manifold extends to a meromorphic map from
the envelope of holomorphy.
\end{Theorem}

\begin{Lemma}
Any meromorphic map $M \to G/H$ to a reductive homogeneous space $G/H$
from a~domain $M$ in Stein manifold extends uniquely
to a meromorphic map $\hat{M} \to G/H$ to its
envelope of holomorphy.
\end{Lemma}
\begin{proof}
By Lemma~\vref{lemma:ClassicalInvariantTheory}, $G/H$ is an
af\/f\/ine algebraic variety. Apply the theorem of Kajiwara and
Sakai (Theorem~\vref{theorem:KajiwaraSakai} above) for
meromorphic functions.
\end{proof}

\subsection{Extension theorems for homogeneous spaces}

The \Thullen extension problem:
\begin{Lemma}[Ivashkovich~\cite{Ivashkovich:2008}]\label{lemma:MeroEqualsHol}
Suppose that $X$ is a complex manifold with locally transitive
biholomorphism group.
 A local biholomorphism to $X$
extends across a subset of \cct to a local
biholomorphism just when it extends across that subset
to a~meromorphic map.
\end{Lemma}
\begin{proof}
 Take a local biholomorphism $f : M \backslash S \to X$
with $M$ a complex manifold and $S \subset M$
a~subset of \cct.
Suppose that $f$ extends meromorphically
to $f : M \to X$ with graph $\Gamma \subset M \times X$.
We only need to extend holomorphically across $S$
locally, so we can assume that $M$ is Stein
and connected and that $X$ is connected.
Therefore $\Gamma$ is connected.

For each holomorphic vector f\/ield $v_X$ on $X$,
def\/ine a holomorphic vector f\/ield
$v_M$, on $M \backslash S$, by setting
$v_M(z)=f'(z)^{-1} v_X(f(z))$ for $z \in M \backslash S$.
Applying Hartogs extension lemma to the coef\/f\/icients of $v_M$
in local coordinates, we can extend $v_M$
uniquely to a holomorphic vector f\/ield~on~$M$.

For each holomorphic vector f\/ield $v_X$ on
$X$, we also def\/ine a holomorphic vector f\/ield
$v_{M \times X}$ on $M \times X$ by
taking $v_{M \times X}=\left(v_X,v_M\right)$.
Let $\mathfrak{g}$ be the Lie algebra
of all holomorphic vector f\/ields on $X$.
This Lie algebra might be inf\/inite dimensional,
and acts transitively on $X$.
The action of $\mathfrak{g}$ on $M \times X$ maps
to the action of $\mathfrak{g}$ on $X$,
and so the orbits in $M \times X$
must be at least as large in dimension
as $X$. Let $\Gamma_0$ be the part of
the graph of $f$ which lies above
$M \backslash S$. Then $H$ is invariant
under this Lie algebra action.
Moreover $\Gamma_0$ is dense in $\Gamma$.
Therefore $\Gamma$ is invariant
under the Lie algebra action.
Each orbit inside $\Gamma$ has
dimension at least that of $X$,
while $\Gamma$ has dimension equal
to that of $X$. The singular locus
of $\Gamma$ is invariant, so is
a union of orbits, each of dimension
equal to that of $\Gamma$.
So $\Gamma$ has empty singular locus, and is a smooth complex
manifold of dimension equal to the dimension
of $X$. Since $\Gamma$ is connected, and all
orbits on $\Gamma$ are open sets, $\Gamma$
is a single orbit.

The projection map $M \times X$ restricted
to $\Gamma$ is a holomorphic surjective map $\Gamma \to M$,
injective on a Zariski open set. If we can show
that $\Gamma \to M$ is a biholomorphism, then we
can invert it to a holomorphic map
$M \to \Gamma$, and then map $\Gamma \to X$
by the other projection, and the composition
$M \to \Gamma \to X$ holomorphically extends $f$.
To prove that $\Gamma \to M$ is a biholomorphism,
we only have to show that it is a local
biholomorphism, since $\Gamma$ is closed
in $M \times X$, so the number of
sheets of $\Gamma \to M$ won't change
at dif\/ferent points of $M$. Since
$\Gamma \to M$ is $\mathfrak{g}$-equivariant,
the set of points at which
$\Gamma \to M$ fails to be a local biholomorphism
must be $\mathfrak{g}$-invariant,
so empty or all of $\Gamma$. However, $\Gamma \to M$
is a local biholomorphism on a dense open set,
so must be a~biholomorphism everywhere.
\end{proof}

The Hartogs extension problem:
\begin{Lemma}[Ivashkovich \cite{Ivashkovich:2008}]\label{lemma:MeroEqualsHolHartogs}
 A local biholomorphism from a domain in a Stein mani\-fold
to a complex homogeneous space
extends to a local biholomorphism
from the envelope of holomorphy just when it extends
to a meromorphic map from the envelope of holomorphy.
\end{Lemma}
\begin{proof}
The proof of Lemma~\vref{lemma:MeroEqualsHol} also works here.
\end{proof}

\begin{table}[b!]
\caption{The complex homogeneous surfaces.}\label{table:HomogeneousComplexSurfaces}\vspace{-6mm}
$$
\begin{array}{lll}
\hline
                                        & \text{symmetry group}  & \text{stabilizer subgroup} \\[-0.8mm] 
\Proj{2}                                & \PSL{3,\C{}}           &
\tableBMatrix{
 a^0_0 & a^0_1 & a^0_2 \\
 0     & a^1_1 & a^1_2 \\
 0     & a^2_1 & a^2_2
}
\\ 
\C{2}                                   & \text{af\/f\/ine}          & \GL{2,\C{}}
\\[-0.5mm] 
\C{2}                                   & \text{special af\/f\/ine}  & \SL{2,\C{}}
\\[-0.5mm] 
\C{2}                                   & \text{translation}     & 0
\\[-0.5mm] 
\C{2}                                   & \text{etc.}            & \text{etc.}
\\[-0.8mm] 
\C{2} \backslash 0                      & \GL{2,\C{}}            &
\tablePMatrix{
1 & b \\
0 & c
}
\\[-0.8mm] 
\C{2} \backslash 0                      & \SL{2,\C{}}            &
\tablePMatrix{
1 & b \\
0 & 1
}
\\[-0.8mm] 
\Proj{1} \times \C{}                    & \PSL{2,\C{}} \times
\left(\C{\times} \rtimes \C{}\right) &
\tableBMatrix{
a & b \\
0 & a^{-1}
}
\times \C{\times}
\\[-0.8mm] 
\Proj{1} \times \C{}                    & \PSL{2,\C{}} \times
\C{} &
\tableBMatrix{
a & b \\
0 & a^{-1}
}
\\[-0.8mm] 
\Proj{1} \times \Proj{1} \backslash \operatorname{diagonal}
                                        & \PSL{2,\C{}}          &
\tableBMatrix{
a & 0 \\
0 & 1/a
}
\\[-0.8mm] 
\Proj{1} \times \Proj{1}                & \PSL{2,\C{}} \times \PSL{2,\C{}}              &
\tableBMatrix{
a & b \\
0 & 1/a
} \times
\tableBMatrix{
c & d \\
0 & 1/c
}
\\ 
\OO{n}, n \ge 1                         & \left(\GL{2,\C{}}/\Z{}_n\right) \rtimes \Sym{n}{\C{2}}^* &
\tablespacer{
\SetSuchThat{
\left(
\tablePMatrix{
a & b \\
0 & d
}, p\right)}{p(1,0)=0}
}
\\[-0.8mm] 
\OO{n}, n \ge 1                         & \left(\SL{2,\C{}}/\left(\pm^n\right)\right) \rtimes \Sym{n}{\C{2}}^* &
\tablespacer{
\SetSuchThat
{
\left(
\begin{pmatrix}
a & b \\
0 & \frac{1}{a}
\end{pmatrix}, p\right)
}{p(1,0)=0}
}\bsep{2ex}
\\ 
\hline
\end{array}
$$\vspace{-6mm}
\end{table}

\begin{Theorem}[Ivashkovich \cite{Ivashkovich:2008}]\label{theorem:KahlerHomogeneous}
Compact K\"ahler homogeneous spaces of dimension at least two
are Hartogs/\Thullen extension targets for local biholomorphisms.
\end{Theorem}
\begin{proof}
Suppose that $X=G/H$ is a compact K\"ahler homogeneous space.
Wang \cite{Wang:1954} proves that such spaces are
complex homogeneous spaces, i.e.\ we can take
$G$  to be a complex Lie group and $H$ a closed
complex Lie subgroup of $G$. We can even assume that
$G$ is the biholomorphism group of $X$. Now apply
Theorems~\ref{theorem:Siu} and~\ref{theorem:Ivashkovich}
and Lemma~\ref{lemma:MeroEqualsHol}.
\end{proof}

\begin{Theorem}
Suppose that $X$ is a complex manifold, and that
the Lie algebra of all holomorphic vector fields on $X$
acts locally transitively. Furthermore suppose that $X$ is
an unbranched covering space of a compact K\"ahler
manifold. Then $X$
is a Hartogs/\Thullen extension target for local biholomorphisms.
\end{Theorem}
\begin{proof}
Again apply Theorems~\ref{theorem:Siu} and~\ref{theorem:Ivashkovich}
and Lemma~\ref{lemma:MeroEqualsHol}.
\end{proof}

 Clearly if $X$ is any Hartogs/\Thullen
target for local biholomorphisms, we can replace $X$ by any complex
manifold with a common covering space, cut out any hypersurface,
replace with a pseudoconvex open set, etc.\ repeatedly and still
have a Hartogs/\Thullen extension target for
local biholomorphisms, so we have a large collection of
complex manifolds to work with. In this paper, we are
only interested in homogeneous Hartogs/\Thullen extension targets.

\begin{Example}\label{example:homogeneousSurfaces}
Let us determine the extension properties of all complex homogeneous surfaces.
Note that in calling a surface \emph{complex homogeneous} we mean
homogeneous under a holomorphic
action of a complex Lie group.
Up to covering, the complex homogeneous surfaces
are presented in Table~\vref{table:HomogeneousComplexSurfaces}.
Huckleberry \cite{Huckleberry:1986}, Mostow \cite{Mostow:1950}
 and Olver \cite[p.~472]{Olver:1995}
provide an introduction to this classif\/ication;
Mostow provides a proof.
The surface $\Proj{1} \times \Proj{1} \backslash \operatorname{diagonal}$ is
the set of pairs of distinct lines in the plane, acted on by linear maps
of the plane.
Consider the usual line bundle $\OO{n} = \OO{1}^{\otimes n} \to \Proj{1}$,
whose f\/ibers are choices of line in $\C{2}$ and homogeneous polynomial
of degree $n$ on that line. In other
words, $\OO{n}$ is the line bundle on $\Proj{1}$
with f\/irst Chern class $n$. We will also denote the total space of
this line bundle by the symbol $\OO{n}$.
This total space is a surface acted on by the group
of linear substitutions of variables. It
is also acted on by the homogeneous
polynomials of degree $n$, by adding
the homogeneous polynomial to the polynomial on any given line.
The group $\Z{}_n$ is the group of scalings of variables by roots
of unity. The group $\pm^n$ is $\pm 1$ if $n$ is even, and $1$ if $n$ is odd.

The symmetry groups listed are all of the
connected complex Lie groups that act transitively
on each given surface, except for the surface $\C{2}$.
On $\C{2}$, there are f\/inite dimensional
Lie groups of all positive dimensions greater
than one, acting transitively. They are
classif\/ied, but the classif\/ication is a
little complicated and irrelevant here; see
Olver \cite[Cases 1.5--1.9, p.~472]{Olver:1995}.
There are also some disconnected
Lie groups acting on these same surfaces
and containing the connected ones listed;
see Huckleberry~\cite{Huckleberry:1986}.

\centerline{\begin{tabular}{lccc}
Homogeneous                      & Hartogs & \Thullen            & local biholomorphism \\
surface                          & target  & target             & \Thullen target \\
$\Proj{2}$                       & \no     & \no                & \yes \\
$\C{2}$                          & \yes    & \yes               & \yes \\
$\C{2} \backslash 0$             & \no     & \no                & \no \\
$\Proj{1} \times \C{}$           & \no     & \no                & \yes \\
$\Proj{1} \times \Proj{1} \backslash \operatorname{diagonal}$
                                 & \yes    & \yes               & \yes \\
$\Proj{1} \times \Proj{1}$       & \no     & \no                & \yes \\
$\OO{n}, n \ge 1$                & \no     & \no                & \no \\
\end{tabular}}

Map $\C{2} \backslash 0 \to \OO{n}$ by taking
a point $z \ne 0$ to the homogeneous polynomial $p_z$
of degree $n$ on the
span of $z$ which takes the value $1$ on $z$. This
map is a local biholomorphism and doesn't extend across the puncture.

To see which complex homogeneous surfaces are Hartogs
or \Thullen extension targets, we start by looking at
covering spaces, and which complex surfaces contain rational curves.
If a~complex surface contains a rational curve, then it isn't a
Hartogs or \Thullen extension target, as we saw in
Example~\vref{example:extensionSubspaces}. Keep in mind that we can
apply Theorem~\vref{theorem:KahlerHomogeneous} to any compact K\"ahler complex
manifold with the same universal covering space as a given space; for example,
$\Proj{1} \times \C{}$ is the universal covering space
of $\Proj{1} \times E$ for any elliptic curve $E$, to which
we can apply Theorem~\ref{theorem:KahlerHomogeneous}. Finally,
notice that
$\Proj{1} \times \Proj{1} \backslash \text{diagonal}$ is an af\/f\/ine variety,
being a reductive homogeneous space:
\[
\Proj{1} \times \Proj{1} \backslash \text{diagonal}
=
 \PSL{2,\C{}}/
\left\{
\begin{bmatrix}
a & 0 \\
0 & \frac{1}{a}
\end{bmatrix}
\right\}.
\]
\end{Example}

\begin{Example}
Now that we have some understanding of local biholomorphisms to
complex homogeneous spaces,
consider immersions to those spaces.
Take two elliptic curves $C_0 = \C{}/\Lambda_0$ and $C_1 = \C{}/\Lambda_1$.
The manifold $M=\Proj{1} \times C_1 \times C_2$ is a complex
homogeneous projective variety. The immersion $\C{2} \backslash 0 \to M$,
$\left(z_1,z_2\right) \mapsto \left(\frac{z_1}{z_2}, z_1 + \Lambda_1, z_2 + \Lambda_2\right)$,
does not extend across the puncture. We can embed $M$ into projective
space to see that immersions to projective space don't always extend.
\end{Example}

\subsection{Extensions of integral maps of invariant dif\/ferential relations}

\begin{Definition}
 Suppose that $G/H$ is a complex homogeneous space.
An \emph{invariant relation} of f\/irst order is a complex
submanifold $R \subset \left(\mathfrak{g}/\mathfrak{h}\right) \otimes \C{n*}$,
for some integer $n$,
invariant under the obvious action of $H \times \GL{n,\C{}}$.
Given any holomorphic map $f : M \to G/H$
from any $n$-dimensional complex manifold $M$,
def\/ine $f^{(1)} : FTM \times_M f^* G \to
\left(\mathfrak{g}/\mathfrak{h}\right) \otimes \C{n*}$
by $f^{(1)}(u,g) = \left(L_g^{-1}\right)'(f(m))f'(m)u^{-1}$.
An \emph{integral map} of $R$ is a holomorphic
map $f : M \to G/H$ of an $n$-dimensional complex
manifold $M$ for which $f^{(1)}$ is valued in $R$.
A relation $R$ will be called a \emph{\Thullen extension relation}
if every integral map of the relation $f : M \backslash S \to G/H$,
with $M$ a complex manifold and $S \subset M$ a subset of \cct,
extends holomorphically across $S$. We
will use the term \emph{Hartogs extension relation} analogously.
\end{Definition}

\begin{Example}
 The set of linear isomorphisms $\C{n*} \to \mathfrak{g}/\mathfrak{h}$
is a Hartogs and a \Thullen extension relation.
\end{Example}

We will see that extension relations
are very closely related to extension problems
for va\-rious geometric structures, giving rise to numerous examples.
Clearly if $R_0 \subset R_1$ are invariant relations and $R_1$
is a \Thullen/Hartogs extension relation, then $R_0$ is too. So naturally
we want to focus on f\/inding maximal Hartogs/\Thullen extension relations.

\section{Extensions of bundles and geometric structures}

\subsection{Extending holomorphic bundles}

The results of this subsection are well known
(see Okonek et al. \cite[Chapter II.1.1]{Okonek/Schneider/Spindler:1980}, Siu \cite{Siu:1974}),
but we provide elementary proofs for completeness.

\begin{Definition}
Suppose that $M$ is a complex manifold, $S \subset M$
a closed subset, $G$ a complex Lie group, and
$E \to M$ a holomorphic principal $G$-bundle.
We will say that $E$ extends
across $S$ to mean that there
is a holomorphic principal $G$-bundle $E' \to M$
and a $G$-equivariant biholomorphism $E \to \left.E'\right|_{M \backslash S}$.
\end{Definition}

\begin{Lemma}\label{lemma:uniqueBundleExtension}
The extension of a holomorphic principal bundle across a
subset of \cct is unique.
\end{Lemma}
\begin{proof}
 Suppose that $E' \to M$ and $E'' \to M$ are two
holomorphic principal $G$-bundles, that $S \subset M$
is a subset of \cct, and that
$F : \left.E'\right|_{M \backslash S} \to \left.E''\right|_{M \backslash S}$
is a~holomorphic principal $G$-bundle isomorphism.
If $F$ extends holomorphically to some open set, then clearly
it does so uniquely. We can therefore replace
$M$ by an open subset on which $E'$ and $E''$ are both holomorphically
trivial. So $F : \left(M \backslash S\right) \times G
\to \left(M \backslash S\right) \times G$
is expressed as $F(z,g)=\left(z,f(z)g\right)$, for some map $f : M \backslash S \to G$.
By Lemma~\vref{lemma:ComplexLieGroups}, complex Lie groups are \Thullen extension
targets, so $f$ extends to a map $f : M \to G$, and this extends
$F$ to an isomorphism $F(z,g)=\left(z,f(z)g\right)$.
\end{proof}

Clearly the same proof as above shows that the Hartogs
extension problem for holomorphic principal
bundles has at most one solution.

\begin{Definition}
Suppose that $M$ is a complex manifold, $S \subset M$
is a closed subset, and
that $E \to M \backslash S$ is a holomorphic f\/iber bundle.
Say that $E$ is \emph{holomorphically trivial along} $S$ if for each $s \in S$,
there is an open set $U \subset M$ containing
$s$ on which $\left.E\right|_{U \backslash S}$ is holomorphically trivial.
\end{Definition}

\begin{Proposition}\label{proposition:fiberBundle}
Any holomorphic fiber bundle
extends across a closed subset if and only if
it is holomorphically trivial along that subset.
\end{Proposition}
\begin{proof}
If such an extension $E'$ exists, then take $U$ any open set containing $s$ on
which $E'$ is trivial.

Conversely, suppose that for every point $s \in S$, we have an open set $U$ containing
$s$ on which $\left.E\right|_{U \backslash S}$ is trivial.
Cover $M \backslash S$ by all possible open sets on which $E$ is trivial
and pick a~trivialization on each.
Let $F$ be a typical f\/iber of $E$.
To each pair $V,W$ of open sets from this cover, associate
a holomorphic transition map $\phi^V_W : \left(V \cap W\right) \times F \to F$.
Then for this particular subset $U$,
def\/ine $\phi^U_V= \phi^{U \backslash S}_V$,
and $\phi^V_U= \phi^V_{U \backslash S}$.
Clearly $\phi^U_V \phi^V_U=\text{id}$
and $\phi^U_V \phi^V_W \phi^W_U= \text{id}$,
so these transition maps determine a unique holomorphic bundle
$E' \to M$ with obvious canonical isomorphism $\left.E'\right|_{M \backslash S}=E$.
\end{proof}

\begin{Corollary}\label{corollary:section}
Any holomorphic principal bundle
extends across a closed nowhere dense subset if and only if
it has a holomorphic section near each point of that subset.
\end{Corollary}

\begin{Lemma}\label{lemma:coveringSpacesExtendCCT}
Covering spaces extend uniquely across subsets of \cct.
\end{Lemma}
\begin{proof}
Suppose that $M$ is a complex
manifold and $S \subset M$ is a subset of \cct.
Pick a covering space $Z \to M \backslash S$.
Pick any point $s \in S$ and any simply connected
open set $U \subset M$ containing $s$.
Then $U \backslash S$ is also simply connected,
so $\left.Z\right|_{U \backslash S} \to U \backslash S$
is holomorphically trivial. Apply
Lemma~\vref{lemma:uniqueBundleExtension} and
Proposition~\vref{proposition:fiberBundle}.
\end{proof}

\begin{Lemma}\label{lemma:DiscreteKernel}
Suppose that $\tilde{G} \to G$
is a Lie group morphism with discrete kernel $K$.
Any holomorphic principal $\tilde{G}$-bundle
$\tilde{E}$ has quotient $E=\tilde{E}/K$.
Moreover, $\tilde{E}$ extends across
subsets of \cct just when $E$ does.
\end{Lemma}
\begin{proof}
Suppose that $M$ is a complex manifold,
that $S \subset M$ is a subset of \cct, and
that $\tilde{E} \to M \backslash S$ is a holomorphic
principal $\tilde{G}$-bundle.
Clearly if $\tilde{E}$ extends over $M$, then
$E=\tilde{E}/K$ does too. Suppose that $E$ extends
over $M$, so has a local section $\sigma$ def\/ined near some
point $s \in S$, say on an open set $U$.
The preimage of $\sigma$ in $\tilde{E}$ is a covering
space of~$U \backslash S$; apply Lemma~\vref{lemma:DiscreteKernel}.
\end{proof}

\begin{Corollary}\label{corollary:ExtendTensoredLineBundles}
Take any two holomorphic line bundles $L_1$ and $L_2$
with a common tensor power,
say $L_1^{\otimes n_1} = L_2^{\otimes n_2}$,
for some integers $n_1$ and $n_2$. Then
$L_1$ extends across a subset of \cct
just when $L_2$ does.
\end{Corollary}
\begin{proof}
The associated principal bundles are both covering spaces
of the associated principal bundle of $L_1^{\otimes n_1}=L_2^{\otimes n_2}$,
so we can apply Lemma~\vref{lemma:DiscreteKernel} to each.
\end{proof}

\begin{Example}
On $\C{2} \backslash 0$, all real rank 2 bundles are smoothly trivial.
There is an inf\/inite dimensional space of
holomorphic line bundles on $\C{2} \backslash 0$ which
are not holomorphically trivial; see Serre \cite[p.~372]{Serre:1966} for proof (without
explicit examples). So smooth category obstructions are not enough
to decide whether holomorphic bundles extend across punctures.
\end{Example}

\subsection{Relative extension problems for holomorphic bundles}

\begin{Definition}
Suppose that $G$ and $G'$ are complex Lie groups
and that $\rho : G \to G'$ is a morphism of complex Lie groups.
Suppose that $E \to M$ is a holomorphic principal right $G$-bundle. Let
$G$ act on $E \times G'$ by having $g \in G$ act
on $\left(e,g'\right)$ to give
\[
r_{g}\left(e,g'\right)=\left(r_g e,  g' \rho\left(g\right) \right).
\]
We denote the quotient as $E'=E \times_G G'$. Let
$G'$ act on $E \times G'$ by
\[
R_{h'} \left(e,g'\right)=\left(e,\left(h'\right)^{-1} g'\right)
\]
for $e \in E$ and $g', h' \in G'$. This
action commutes with the $G$-action, so
descends to a right $G'$-action on $E'=E \times_{G} G'$.
Moreover $E' \to M$ is a principal right holomorphic $G'$-bundle.

Consider the map $e \in E \mapsto \left(e,1\right) \in E \times G'$.
Compose this map with the obvious quotient
map $E \times G' \to E \times_G G'=E'$
to make a $G$-equivariant
bundle map $p : E \to E'$.
We def\/ine $q : E \times G' \to G'$ by
$q\left(e,g'\right)=g'$.
The map $q$ is $G$-equivariant, so descends to a
map $q : E' \to G'/\rho\left(G\right)$.
Moreover, the composition $qp$
is the constant map to the $\rho\left(G\right)$ coset. Conversely,
the $\rho\left(G\right)$-subbundle
$E/\ker \rho = p\left(E\right) \subset E'$ is precisely
$E/\ker \rho = q^{-1} \rho\left(G\right)$.
\end{Definition}

\begin{Lemma}\label{lemma:Relative}
Take a morphism $\rho : G \to G'$ of complex Lie
groups.
Suppose that $E$ is a~principal $G$-bundle.
Let $E' = E \times_G G'$ as above.
If $E$ extends across a subset
then $E'$ extends as $E'=E \times_G G'$.
If the kernel of $\rho$ is discrete and
the image of $\rho$ is a closed subgroup of $G'$,
and if holomorphic maps to $G'/\rho(G)$
extend across subsets of \cct, then
$E$ extends across a subset of \cct just when $E'$ does.
\end{Lemma}
\begin{proof}
Suppose that $M$ is a complex manifold
and $S \subset M$ is a subset of \cct.
Suppose that $E' \to M \backslash S$ extends
holomorphically to a principal $G'$-bundle
$E' \to M$, and that the image of $\rho$ is a closed subgroup of $G'$,
and that $G'/\rho\left(G\right)$ is a \Thullen extension target.
We can replace $M$ with a small
neighborhood of a point $s \in S$, in which
$E' \to M$ is trivial: $E' = M \times G'$.

Denote a typical point of $E'$ as
$\left(z,g'\right)$ with $z \in M$ and $g' \in G'$.
Write $q$ as
$q \left(z,g'\right) \in G'/\rho\left(G\right)$.
Because $G'/\rho\left(G\right)$ is a
\Thullen extension target, we can extend the map
$q$ to $E'$ by extending it to be
holomorphic for $z \in M$ for each f\/ixed $g'$.
Then $q^{-1} \rho\left(G\right)$
is a holomorphic $G/\ker \rho$-subbundle of $E'$ extending
$E/\ker \rho$. By Lemma~\vref{lemma:DiscreteKernel},
$E$ extends across $S$.
\end{proof}

\subsection{Extending bundles via a connection}

\begin{Example}
Af\/f\/ine connections are given in local coordinates
by Christof\/fel symbols, which are holomorphic
functions. The Christof\/fel symbols extend holomorphically
across subsets of \cct by Hartogs lemma. Therefore
af\/f\/ine connections holomorphically extend across subsets of \cct.
We have already generalized this result to holomorphic
connections on bundles
in Lemma~\vref{lemma:holomorphicConnectionsExtend}.
\end{Example}

\begin{Example}
If $G \subset \GL{n,\C{}}$ is a closed complex Lie subgroup,
and $E \to M \backslash S$ is a $G$-structure equipped with a torsion-free
connection, then the connection extends to a torsion-free
connection of $FTM=E \times_G \GL{n,\C{}}$ with holonomy
in $G$. The connection extends over $S$ by the last
example. The holonomy of a loop $\gamma$ passing through $S$
is the limit of the holonomy of any sequence of loops
$\gamma_j$ converging to $\gamma$ uniformly. We can
choose each $\gamma_j$ to avoid $S$,
so the holonomy of $\gamma_j$ lies in $G$.
Therefore the parallel transport on $FTM$ preserves
a foliation by $G$-subbundles, one leaf $E'$ of which contains and
therefore extends $E$ from $E \to M \backslash S$ to $E' \to M$.
\end{Example}

\begin{Proposition}[Buchdahl and Harris %
\cite{Buchdahl/Harris:1999}]\label{proposition:connection}
A holomorphic principal bundle or vector bundle on a complex
manifold extends across a subset of \cct just when
it admits a holomorphic connection near each point
of that subset.
\end{Proposition}
\begin{Remark}
This solves the \Thullen extension problem for
holomorphic bundles with holomorphic connections;
the Hartogs extension problem is unsolved.
\end{Remark}
\begin{proof}
Suppose that $E \to M \backslash S$ is a holomorphic principal bundle,
with $M$ a complex manifold and $S \subset M$
a subset of \cct.
We can assume that $M$ is a ball $B \subset \C{n}$,
and that the holomorphic connection
is def\/ined on all of $E$. Pick any point $z_0 \in B \backslash S$.
Each complex line through $z_0$ intersects $B$ in a~disk. Each complex line also intersects $B \backslash S$ in
a~disk, except for those lines which intersect points
of $S$, which yield disks with f\/initely many punctures.
Parallel transport between any two
points of a disk is well def\/ined, because the connection
is holomorphic. There could be monodromy around a punctured
disk, but since the punctured disk is a limit
of unpunctured disks (coming from the nearby complex lines
through~$z_0$), the monodromy is trivial.
Therefore if we pick an initial point of the f\/iber
$E_{z_0}$, we can parallel transport it around all of the
disks through $z_0$, obtaining a global holomorphic
section of $E$ over~$B \backslash S$.
By Corollary~\vref{corollary:section}, the bundle is holomorphically
trivial. Extend the connection across $S$
by writing it out in local coordinates and applying Hartog's
extension lemma to the Christof\/fel symbols.

For a vector bundle, consider the associated
principal bundle.
\end{proof}

\begin{Definition}
 Suppose that $M$ is a complex manifold, $S \subset M$
is an analytic subset,
$G$ is a~complex Lie group,
and $E \to M \backslash S$ is a holomorphic
principal $G$-bundle. Each open set $U \subset M$ containing
$S$ has an Atiyah class
\[
a\left(\left.E\right|_{U \backslash S}\right)
\in
\Cohom{1}{U \backslash S, T^*\left(U \backslash s\right)
\otimes
\left(
\left.E\right|_{U \backslash S} \times_G \mathfrak{g}
\right)
},
\]
and these Atiyah classes pullback under inclusions of open sets.
For each point $s \in S$, def\/ine the
\emph{Atiyah class} of $E$ at $s$, written
$a(E,s)$, to be the inverse limit of Atiyah
classes $a\big(\left.E\right|_{U \backslash S}\big)$
over all open sets $U \subset M$ containing $s$.
\[
a\left(E,s\right)
\in
\ilim_{s \in U}
\Cohom{1}{U \backslash S, T^*\left(U \backslash s\right)
\otimes
\left(
\left.E\right|_{U \backslash S} \times_G \mathfrak{g}
\right)
}.
\]
If $S$ is subset of \cct
then $a(E,s)=0$ for all $s \in S$ just when $E$ extends across $S$.
\end{Definition}

\subsection{Extending geometric structures by extending bundles}

\begin{Definition}
We will say that $G/H$-geometries \emph{extend across
subsets of \cct} to mean that if $M$ is a complex manifold,
$S \subset M$ a subset of \cct,
and $E \to M \backslash S$ is a
holomorphic $G/H$-geometry
then $E$ extends to a~unique holomorphic $G/H$-geometry
on $M$. Equivalently, we
will refer to solving the \Thullen extension
problem for $G/H$-geometries.
\end{Definition}
Similar terminology will be used for f\/irst order structures, etc.

\begin{Theorem}\label{theorem:CartanGeometry}
The underlying holomorphic principal bundle of a holomorphic Cartan geomet\-ry
extends across a subset of \cct to a holomorphic principal bundle
just when the Cartan geometry extends across the subset.
\end{Theorem}
\begin{proof}
The result is local, so we can assume that $M$ is a ball $B \subset \C{n}$,
and that $E \to B \backslash S$ is a Cartan geometry, extending as a
principal bundle to some  $E' \to B$, which we can assume is
holomorphically trivial, $E'=B \times H$.
Denote points of $E'$ as $\left(z,h\right)$, $z \in B$ and $h \in H$.
Take linear coordinates $z^1, z^2, \dots, z^n$ on $\C{n}$.
From part (3) of the def\/inition of a Cartan geometry
(Def\/inition~\vref{def:CartanConnection}),
$\omega$ restricts to $h^{-1} \, dh$ on each f\/iber $\left\{z\right\} \times H$.
Therefore $\omega - h^{-1} \, dh$ is a~multiple of $dz^1, dz^2, \dots, dz^n$,
say $\omega = h^{-1} \, dh + \Ad(h)^{-1} \left( \Gamma_j(z,h) \, dz^j \right)$,
with $\Gamma_j(z,h)$ holomorphic for $z$ away from $S$ and valued in $\mathfrak{g}$.
From part (1) of the def\/inition of a Cartan geometry, $\Gamma_j(z,h)$
is independent of $h$, say $\Gamma_j(z)$, holomorphic for $z \ne 0$. Therefore
$\Gamma_j(z)$ admits a unique extension to a holomophic function on the
ball $B$, extending $\omega$ to a holomorphic 1-form on $B \times H=E'$. The
properties (1) and (3) of Cartan connections then follow immediately.
In terms of a basis $\left\{e_i\right\}$ of $\mathfrak{g}$,
we can write $\omega = \omega^i e_i$.
The tangent bundle of $B \times H$ is holomorphically trivial,
so $B \times H$ admits a global holomorphic volume form.
Pick any holomorphic volume form $\Omega$ on
$B \times H$. The holomorphic function
\[
 f =
\frac{\Omega}%
{\omega^1 \wedge \omega^2 \wedge \dots \wedge \omega^{\dim G}}
: B \backslash S \to \C{}
\]
extends across $S$ by Hartogs lemma, so $\omega$ satisf\/ies
property (2) of a Cartan connection.
\end{proof}

\begin{Theorem}\label{theorem:CartanGeometryHartogs}
Consider a holomorphic Cartan geometry
defined on a domain in Stein mani\-fold.
The underlying holomorphic principal bundle of the Cartan geometry
extends to a holomorphic principal bundle on the envelope
of holomorphy just when the Cartan geometry extends.
\end{Theorem}
\begin{proof}
Suppose that $\pi : E \to M$ is a holomorphic Cartan
geometry modelled on $G/H$ with Cartan connection $\omega$.
Suppose that $M$ is a domain in a Stein manifold
with envelope of holomorphy $\hat{M}$.
Suppose that $E$ extends to a holomorphic
principal $H$-bundle $E' \to \hat{M}$.
The Cartan connection $\omega$ is a holomorphic
section of the vector bundle $T^*E \otimes \mathfrak{g}$ on $E$.
Being $H$-equivariant, $\omega$ is also a
section of the vector bundle
$\left(T^*E \otimes \mathfrak{g}\right)^H \to M$,
whose local sections are $H$-equivariant local
sections of $T^*E \otimes \mathfrak{g}$. This
latter vector bundle extends to
the holomorphic vector bundle
$\left(T^*E' \otimes \mathfrak{g}\right)^H \to \hat{M}$.
By Proposition~\vref{proposition:VectorBundleSections},
the section $\omega$ extends to a section of this
vector bundle, which we will also denote as $\omega$.
This section is thus an $H$-invariant
1-form valued in $\mathfrak{g}$, the f\/irst property
of a Cartan connection (see Def\/inition~\vref{def:CartanConnection}).

Let $Z$ be the set of points $e \in E'$ at
which $\omega_e : T_e E' \to \mathfrak{g}$
is \emph{not} a linear isomorphism. Clearly~$Z$ is a hypersurface, given by the one equation
$\det \omega_e=0$. Moreover, this hypersurface
is $H$-invariant, so projects to a hypersurface
in $\hat{M}$. This hypersurface doesn't intersect
$M$, so is empty by Lemma~\vref{lemma:hypersurfaceIntersections}.
Therefore $\omega$ satisf\/ies the second property
of a Cartan connection.

Over $M$, $\omega$ satisf\/ies $\vec{A} \hook \omega = A$,
for any $A \in \mathfrak{h}$. By analytic continuation,
this must also hold over $\hat{M}$, the third and f\/inal property
of a Cartan connection.
\end{proof}

\begin{Theorem}\label{theorem:GeometricStructure}
The underlying holomorphic principal bundle of a
holomorphic first order structure
extends across a subset of \cct
to a holomorphic principal bundle
just when the first order structure extends across the subset.
\end{Theorem}
\begin{proof}
The problem is local, so we can assume that $M$ is a ball $B \subset \C{n}$,
and that $E' \to M$ is holomorphically trivial, $E'=B \times G$.
Denote points of $E'$ as $\left(z,g\right)$, $z \in B$ and $g \in G$.
Take linear coordinates $z^1, z^2, \dots, z^n$ on $\C{n}$.
We have a map
$\phi : E \to FT\left(M \backslash S\right),
\phi\left(z,g\right)=\left(z,u\left(z,g\right)\right)$,
def\/ined for $z \in B \backslash S$,
with $u(z,g) \in \GL{n,\C{}}$. Moreover, $u(z,g)=g^{-1}u(z,1)$. So
we can consider the map $u(z)=u(z,1)$ as a map $u : B \backslash S \to \GL{n,\C{}}$.
Clearly $u$ extends to a matrix-valued function on $B$.
Moreover, so does $u^{-1}$, so clearly $u$ extends to a
map $u : B \to \GL{n,\C{}}$, and we extend $\phi$ to $E'$ by
$\phi(z,g)=\left(z,g^{-1}u(z)\right)$.
\end{proof}

\begin{Conjecture}
The underlying holomorphic principal bundle of a
holomorphic first order structure with trivial kernel
extends from a domain $M$ in a Stein manifold
to a holomorphic principal bundle
on the envelope of holomorphy of $M$
just when the first order structure extends as well.
\end{Conjecture}

\section{Extending f\/irst order structures}

\subsection{Inextensible examples}

\begin{Example}\label{example:proj}
Let $G \subset \GL{n,\C{}}$ be
the stabilizer of a nonzero vector $v_0 \in \C{n}$.
A $G$-structure on a manifold $M$
is precisely a nowhere vanishing
vector f\/ield on $M$. Suppose that $\phi : E \to FTM$
is a $G$-structure on $M$, with underlying
principal right $G$-bundle $\pi : E \to M$.
Take any point $e \in E$.
Since $G \subset \GL{n,\C{}}$ is a
subgroup, and $\phi$ is a $G$-bundle
morphism, this point $e$ is identif\/ied with a
point $u=\phi(e) \in FTM$.
This point $u$ is a linear
isomorphism $u : T_m M \to \C{n}$, where
$m=\pi(e)$.
Consider the vector $X=u^{-1}\left(v_0\right) \in T_m M$.
Denote this vector $X(e)$.
Clearly since $v_0 \ne 0$ and $u$ is a linear
isomorphism, $X(e) \ne 0$.

Under $G$-action, $\phi(eg)=g^{-1}\phi(e)$, so
$X(eg)=u^{-1} g^{-1}\left(v_0\right)=u^{-1}\left(v_0\right)=X(e)$.
Therefore $X(e)$ depends only on the point $m \in M$: $X=X(m)$,
i.e. $X$ is a vector f\/ield on $M$. By local
triviality of $\pi : E \to M$, $X$ is a holomorphic
vector f\/ield, nowhere vanishing.

Conversely, suppose that $X$ is a nowhere
vanishing vector f\/ield on a manifold $M$.
Let $E$ be the set of all pairs
$e=(m,u)$ for which $m \in M$ and
$u : T_m M \to \C{n}$ is a linear
isomorphism satisfying
$u(X(m))=v_0$. Clearly $E \subset FTM$
is a $G$-structure. So we have an
isomorphism between the category
of $G$-structures on $M$ and the category
of nowhere vanishing vector f\/ields
on $M$.

Vector f\/ields extend across
subsets of \cct, by Hartogs extension
lemma applied in local coordinates to the component functions
of the vector f\/ield. (They even extend to the envelope
of holomorphy, if there is one.)
Clearly the associated $G$-structure extends across such
a subset just when the vector f\/ield
extends without zeroes.
As an example, take
an invertible $n \times n$ matrix $A$,
and let $X$ denote the vector f\/ield $X(z) = Az$ on
$\C{n}$. Clearly $X$ doesn't vanish
except at $z=0$. Moreover, there is a
unique holomorphic extension of $X$
from $\C{n} \backslash 0$ to $\C{n}$.
Let $S=\left\{0\right\}$, $M=\C{n}$, and let $E$ be the
$G$-structure on $M \backslash S$ associated
to $X$. Then $E$ does not extend
holomorphically as a $G$-structure
to $M$. By Theorem~\vref{theorem:GeometricStructure},
the holomorphic bundle $E \to \C{n} \backslash 0$
does not extend to a holomorphic principal
bundle over  $\C{n}$. So we cannot always
solve the \Thullen or Hartogs extension problems
for f\/irst order structures.
\end{Example}

\begin{Example}\label{example:hyperplanes}
Consider the 1-form $\alpha=\sum z^i \, dz^i$ on $\C{n}$.
On $\C{n} \backslash 0$, take the
hyperplane f\/ield $\alpha=0$. Let $G$ be
the group of linear maps on $\C{n}$ preserving
a hyperplane, say $\C{n-1} \subset \C{n}$.
Let $E$ be the set of pairs
$\left(z,u\right)$ with $z \in M \backslash S$
and $u \in \GL{n,\C{}}$ for which
$u$ takes the hyperplane $\left(\alpha=0\right)$ to the
f\/ixed hyperplane $\C{n-1}$. Clearly $E$ is
a $G$-structure on $M \backslash S$.

Let us show that this $G$-structure does not extend
across the puncture at $0$. If this $G$-structure
extends across the puncture, then we can take a local
section, say $u(z)$, def\/ined near~$0$, and def\/ine
a 1-form $\beta$ by $v \hook \beta=dz^1(u(v))$.
This 1-form $\beta$ doesn't vanish at any point,
and annihilates tangent vectors
precisely on the hyperplanes, as does $\alpha$, so
$\beta = h   \alpha$ for some nonzero function~$h$
away from $0$. Extend $\beta$ and $h$ and $1/h$
to $0$ by the Hartogs extension lemma.
So~$\beta$ vanishes at~$0$, a
contradiction. So we cannot always solve the
Hartogs or \Thullen extension problems for
hyperplane f\/ields.
\end{Example}

\begin{Example}
Generalizing the previous example, we can take
any closed complex subgroup $G \subset \GL{n,\C{}}$,
any nonconstant map $f : \Proj{1} \to \GL{n,\C{}}/G$
(assuming there is a nonconstant map, which is a complicated
constraint on the choice of $G$).
For a point $z \in \C{2} \backslash 0$,
write $[z]$ for the complex line through $z$ and
$0$, mapping $\C{2} \backslash 0 \to \Proj{1}$.
We map $F : \left(\C{2} \backslash 0\right) \times \C{n-2} \to \GL{n,\C{}}/G$
by $F(z,w)=f([z])$.
Take the bundle $\GL{n,\C{}} \to \GL{n,\C{}}/G$ and
let $E=F^* \GL{n,\C{}}$ be the pullback.
The $G$-structure can't extend
across the puncture at $z=0$ because the function $f$ on $E$ is def\/ined
on $\C{2} \backslash 0$ and is constant
along all lines through $z=0$, with dif\/ferent constants
along dif\/ferent lines. Intuitively, this tells us that
we cannot solve the Hartogs or \Thullen extension
problems for $G$-structures for ``large'' subgroups $G \subset \GL{n,\C{}}$
(for example, parabolic subgroups). For such ``large'' subgroups,
we will need to use additional hypotheses on torsion
(see Section~\vref{section:FullTorsion}).
\end{Example}

\begin{Example}
For some choices of group $G$, all $G$-structures extend
across subsets of \cct. For example, an $\SO{n,\C{}}$-structure
(a.k.a. a holomorphic Riemannian metric)
is given in local coordinates by a symmmetric matrix $g=\left(g_{ij}\right)$
of holomorphic functions with $\det g \ne 0$.
The functions will extend holomorphically over such subsets
as will the function $1/\det g$, by
Hartogs lemma. Therefore any holomorphic
Riemannian metric extends holomorphically over subsets of \cct.
Exactly the same trick works for holomorphic symplectic structures.
Similar remarks apply to the Hartogs extension problem for
these structures.
\end{Example}

\begin{Example}
 On $\C{4}$, with coordinates $z^0$, $z^1$, $w^0$, $w^1$, let $\omega=dz^0 \wedge dz^1 + dw^0 \wedge dw^1$
and let $L(z,w)=(iz,-iw)$. Then the span of a pair of vectors
\[
 (z,w), L(z,w)
\]
is a Lagrangian 2-plane for $\omega$, unless
$z=0$ or $w=0$, in which case it is a sub-Lagrangian line.
On $\C{4} \backslash \left( \left(z=0\right) \cup \left(w=0\right) \right)$,
this Lagrangian foliation is holomorphic, and does not extend across
$z=0$ or $w=0$. Consider the usual map $\pi : z \in \C{4} \backslash 0 \to \C{\times}z \in \Proj{3}$.
Let $\alpha = \left(z,w\right) \hook \omega$.
The hyperplane f\/ield $V_{(z,w)} = \ker \alpha$ has the
f\/ibers of $\pi$ as Cauchy characteristics, and so descends to a hyperplane
f\/ield on $\Proj{3}$, which is a holomorphic contact structure. The
Lagrangian 2-planes project to Legendre lines, foliating $\Proj{3}$ away
from the two lines $[z=0]$ and $[w=0]$. This Legendre foliation does not
extend across those two lines.
\end{Example}

\begin{Example}\label{example:PlaneFields}
A \emph{plane field} on a manifold $M$ is a vector subbundle
$V \subset TM$.
An \emph{Engel} 2-plane f\/ield $V \subset TM$ on a 4-fold $M$
is a 2-plane f\/ield so that near each point $m \in M$,
there are local sections $X$ and $Y$ for which
$X, Y, [X,Y]$ and $[[X,Y],Y]$ are linearly
independent.
Every Engel 2-plane f\/ield is locally isomorphic to every other;
see
Ehlers et al. \cite{Ehlers/Koiller/Montgomery/Rios:2005},
Kazarian, Montgomery and Shapiro \cite{Kazarian/Montgomery/Shapiro:1997},
Vogel \cite{Vogel:2006},
or
Zhitomirski{\u\i} \cite{Zhitomirskii:1990}.
In particular, all Engel
2-plane f\/ields are locally isomorphic to
the plane f\/ield
\[
dy - p \, dx=0, \qquad dp - q \,dx =0.
\]
In other words, this is the plane f\/ield
spanned by the vector f\/ields
\[
\pd{}{q}, \qquad \pd{}{x} + p \, \pd{}{y} + q \, \pd{}{p}.
\]
However, it is easy to check that the plane f\/ield
\[
dy - pq \, dx=0, \qquad dp - q^2 \, dx=0,
\]
spanned by the vector f\/ields
\[
\pd{}{q}, \qquad \pd{}{x} + pq \, \pd{}{y} + q^2 \, \pd{}{p}
\]
is also an Engel 2-plane f\/ield except on the
subset $p=q=0$.
On that subset, the 2-plane f\/ield fails to be Engel,
as the required brackets of any local sections
lose their linear independence. Therefore Engel 2-plane f\/ields do
\emph{not} always extend holomorphically
across complex codimension~2 subsets
to Engel 2-plane f\/ields.
\end{Example}

\begin{Conjecture}
Holomorphic Engel $2$-plane fields extend holomorphically
across subsets of complex codimension~$2$
or more to holomorphic $2$-plane fields,
and across subsets of complex codimension~$3$
or more to holomorphic Engel $2$-plane fields.
\end{Conjecture}

\subsection{Contact structures}\label{subsection:ContactStructures}

\begin{Definition}
A \emph{hyperplane field} on a complex manifold $M$ is a holomorphic line
subbundle of the holomorphic cotangent bundle. If $L \subset T^*M$
is a hyperplane f\/ield, then locally $L$ is spanned by a nonzero
1-form, say $\alpha$.
The hyperplane f\/ield is a \emph{contact structure}
if $\alpha \wedge \left(d\alpha\right)^n \ne 0$, where~$M$ has dimension $2n+1$.
\end{Definition}
 A contact structure can be considered a f\/irst order structure with torsion
condition; see Example~\vref{example:hyperplanes}.

\begin{Theorem}\label{theorem:ContactStructuresExtend}
 Holomorphic contact structures extend holomorphically across
subsets of \cct.
\end{Theorem}

To prove this theorem, we will need to prove two minor results.

\begin{Lemma}\label{lemma:contactForm}
Suppose that $M$ is a complex manifold, that $S \subset M$ is a subset of \cct,
and that $M \backslash S$ bears a contact structure.
If, for each point $s \in S$, the contact structure
is spanned by a nonzero $1$-form $\alpha$ defined
in an open set of the form $U \backslash S$, where
$U$ is an open subset of $M$ containing $s$, then
the contact structure extends holomorphically
to a unique contact structure on $M$.
\end{Lemma}

\begin{proof}
Imagine a holomorphic contact structure on $M^{2n+1} \backslash S$
with a choice of holomorphic contact form in a neighborhood of $s$,
say $\alpha$, so that $\alpha \wedge d \alpha^n \ne 0$.
Applying Hartogs extension to the coef\/f\/icients of $\alpha$
in local coordinates near $s$, the contact
form $\alpha$ extends uniquely as a~holomorphic 1-form
across $S \cap U$.
Take $\Omega$ a holomorphic volume form def\/ined near $s$,
and let $f=\Omega/\alpha \wedge d \alpha^n$. By Hartog's
extension lemma, $f$ extends to a holomorphic function
across $S \cap U$, so $\alpha$ extends to a contact form.
\end{proof}

\begin{Corollary}\label{corollary:ContactBundle}
Suppose that $M$ is a complex manifold, that $S \subset M$ is
a subset of \cct,
and that $M \backslash S$ bears a contact structure
$L \subset T^*\left(M \backslash S\right)$.
Then the line bundle $L$ extends
across $S$ as a holomorphic line bundle
if and only if the contact
structure extends across
$S$ as a holomorphic contact structure.
\end{Corollary}
\begin{proof}
If the line bundle $L$ extends as a holomorphic line bundle
across the puncture, then $L$ is locally trivial,
so we can choose a local section, $\alpha$,
and apply Lemma~\vref{lemma:contactForm}.
\end{proof}

\begin{Example}
On the other hand, as we saw in Example~\vref{example:hyperplanes},
the hyperplane f\/ield $z_i \, dz_i=0$
does not extend across $0$ as a hyperplane f\/ield.
Its associated line bundle is trivial over
$\C{n} \backslash 0$, having global section $z_i \, dz_i$.
Therefore the line bundle extends holomorphically
across $0$. This hyperplane f\/ield is not a contact plane f\/ield.
Clearly the contact condition simplif\/ies
the extension problem. More generally, we should expect
nonvanishing torsion of a f\/irst order structure
to be helpful in extension problems.
\end{Example}

Now we return to proving Theorem~\ref{theorem:ContactStructuresExtend}.

\begin{proof}
Suppose that $S \subset M$ is a subset of \cct in a complex manifold,
and that $M \backslash S$ bears a contact structure
$L \subset T^*\left(M \backslash S\right)$.
Consider the inclusion $\iota : L \to T^*\left(M \backslash S\right)$
as a linear map. Take the transpose $\alpha = \iota^t :
T\left(M \backslash S\right) \to L^{\otimes -1}$,
which is a 1-form valued in $L^{\otimes -1}$.
The expression $\alpha \wedge \left(d \alpha\right)^n$
is a section of the bundle
$K_{M \backslash S} \otimes L^{\otimes(-n-1)}$,
nowhere vanishing, so an isomorphism $L^{\otimes(n+1)} \to K_{M \backslash S}$.
Clearly $K_{M \backslash S}$ extends to~$K_M$. Therefore $L^{\otimes(n+1)}$ extends
across $S$. Apply Corollary~\vref{corollary:ExtendTensoredLineBundles} to
conclude that~$L$ extends holomorphically across $S$.
Apply Corollary~\vref{corollary:ContactBundle} to conclude
that the contact structure extends.
\end{proof}

\begin{Example}
Consider the real analytic
contact structure
\[
\cos z \, dx - \sin z \, dy = 0
\]
on $\R{3}$, and then
compactify $\R{3}$ to the 3-sphere. Near
the point at inf\/inity, the contact planes
wind inf\/initely often around certain great
circles, so the contact structure doesn't
extend. Perhaps there is
no complex analytic analogue of
overtwisted contact structures.
\end{Example}

\begin{Proposition}\label{proposition:HartogsSemiExtend}
For every holomorphic subbundle $V \subset TM$
of the tangent bundle of a~domain $M$
in a Stein manifold, there is
a subset of \cct $S \subset \hat{M}$
of the envelope of holomorphy $\hat{M}$
of $M$, so that $V$ extends
uniquely to a holomorphic subbundle of
$T\big(\hat{M} \setminus S\big)$.
\end{Proposition}
\begin{Remark}
The author would like to thank Sergei
Ivashkovich for providing this proof.
\end{Remark}
\begin{proof}
Suppose that $V \subset TM$
is a holomorphic subbundle,
say of rank $k$,
and $M$ is a domain in a Stein
manifold, and $\hat{M}$ is the
envelope of holomorphy of $M$.
There is a holomorphic
embedding $\hat{M} \to \C{N}$
for some integer $N$;
see H\"ormander~\cite[Theorem 5.3.9, p.~135]{Hormander:1990}.

Let $f : M \to \Gr{k}{N}, m \mapsto V_m$.
By Theorem~\vref{theorem:Ivashkovich},
$f$ extends to a meromorphic map
$f : \hat{M} \to \Gr{k}{N}$. By analytic
continuation, the equation $f(m) \subset T_m M$
for $m \in M$ ensures that $f(m) \subset T_m \hat{M}$
for $m \in \hat{M}$.
The indeterminacy locus of
the map $f$ is of \cct.
\end{proof}

An alternative proof making use of more
elementary results, on meromorphic functions
rather than maps, could proceed by
taking af\/f\/ine coordinates on the
Grassmannian and then applying
the Kajiwara--Sakai theorem
(Theorem~\vref{theorem:KajiwaraSakai}).

\begin{Theorem}
Holomorphic contact structures extend uniquely
from any domain $M$ in Stein manifold to the
ennvelope of holomorphy of $M$.
\end{Theorem}
\begin{Remark}
The author would like to thank
Sergei Ivashkovich for providing this proof.
\end{Remark}
\begin{proof}
By Proposition~\ref{proposition:HartogsSemiExtend},
the problem reduces to extension
in complex codimension 2, so the result
follows from
Theorem~\vref{theorem:ContactStructuresExtend}.
\end{proof}

\begin{Example}[Robert Bryant]
Let $X$ be a complex a 2-torus.
Naturally $\Proj{}T^* X$ has the obvious
contact structure. But $\Proj{}T^* X$
also has various holomorphic
2-plane f\/ields.
Note that $\Proj{}T^*X=X \times \Proj{1}$.
In linear coordinates $z$, $w$ on $X$,
and af\/f\/ine coordinate $p$ on $\Proj{1}$,
pick a rational function $f(p)$
and consider the 2-plane f\/ield
$dw - f(p) \, dz=0$.
Note that near points where $f(p)=\infty$,
we can write this 2-plane f\/ield
as
\[
\frac{dw}{f(p)} - dz = 0,
\]
so the 2-plane f\/ield is holomorphic
at all points of $X \times \Proj{1}$.
The 2-plane f\/ield is invariant
under translations of $X$.

The critical points of $f$ are
precisely the points where
the 2-plane f\/ields fails to
be a contact structure.
These points form a union of hypersurfaces,
each one a torus. These
tori are homogeneous under
the translation action of $X$. So we
have many examples of holomorphic
2-plane f\/ields
on a homogeneous smooth projective
variety, which are contact
structures except on some
disjoint subvarieties,
each of which is homogeneous.
The pseudogroup of local isomorphisms
of any of these 2-plane f\/ields
will active transitively on the
dense open set where the 2-plane
f\/ield is a contact structure,
and also on the various smooth
hypersurfaces on which the
2-plane f\/ield fails to be a contact
structure.
\end{Example}

\subsection{Reducing to a homogeneous space extension problem}

\begin{Definition}
Take a complex Lie group $G$ and complex representation $\rho : G \to \GL{V_0}$
with closed image. Suppose that $\phi : E \to FTM$ is a holomorphic $G$-structure,
with $FTM$ the $V_0$-valued frame bundle.
We def\/ine $\Gromov{\phi} : FTM \times_M E \to \GL{V_0}$
by $\Gromov{\phi}(u,e)=u\phi(e)^{-1}$.
Under $G$-action, $\Gromov{\phi}\left(u,r_g e\right)=\Gromov{\phi}(u,e)g$.
Therefore we can quotient by right $G$-action, to produce
a map $\Gromov{\phi} : FTM \to \GL{V_0}/\rho(G)$, with
$\Gromov{\phi}(u)=u\phi(e)^{-1}\rho(G)$ for some element
$u \in FT_m M$ and $e \in E_m$. If we change the choice
of $e \in E_m$, say to $eg$,
we change $\phi(e)^{-1}$
to $\phi\left(eg\right)^{-1}=\phi\left(e\right)^{-1} \rho(g)$,
not af\/fecting $\Gromov{\phi}(u)$.
So $\phi : E \to FTM$ determines $\Gromov{\phi} : FTM \to \GL{V_0}/\rho(G)$,
equivariant under right $\GL{V_0}$-action.
Moreover, the composition $\Gromov{\phi} \phi$
is the constant map to the $\rho(G)$ coset. The $\rho(G)$-structure
$\phi(E) \subset FTM$ (i.e.\ the underlying embedded f\/irst order
structure) is precisely $\Gromov{\phi}^{-1} \rho(G)$.
\end{Definition}

\begin{Theorem}\label{theorem:ReductiveGstructuresHartogsExtension}
Suppose that $G \subset \GL{n,\C{}}$ is a reductive
algebraic group.
Suppose that $M$ is a domain in a Stein manifold.
Then every holomorphic $G$-structure on $M$ extends
uniquely to a holomorphic $G$-structure on the
envelope of holomorphy of $M$.
\end{Theorem}
\begin{proof}
A $G$-structure on $M$ is equivalent to a section of
$FTM/G \subset FT\hat{M}/G$.
The total space of $FT\hat{M}/G$ is
a Hartogs extension target by
Corollary~\vref{corollary:ExtendReductiveReduction}.
Therefore if $s : M \to FTM/G$ is a $G$-structure,
then $s$ extends uniquely to a holomorphic map
$s : \hat{M} \to FT\hat{M}/G$.
Let $\pi : FT\hat{M}/G \to \hat{M}$ denote
the bundle map. Then $\pi s$ is the identity
on $M$, and therefore by analytic continuation
is the identity on $\hat{M}$. So the
extension is also a holomorphic $G$-structure
\end{proof}

\begin{Theorem}\label{theorem:ReductiveGstructures}
Suppose that $\rho : G \to \GL{n,\C{}}$ is a representation
with closed image and discrete kernel.
Then all holomorphic $G$-structures extend across
subsets of \cct if and only if
all holomorphic maps to $\GL{n,\C{}}/\rho(G)$
from $n$-folds extend across
such subsets.
\end{Theorem}
\begin{proof}
We can quotient by the kernel, by Lemma~\vref{lemma:Relative},
so assume that $G \subset \GL{n,\C{}}$
is a closed subgroup.
Suppose that all holomorphic maps from
$n$-folds to $\GL{n,\C{}}/G$ extend
across subsets of \cct.
We need only prove the result locally.
The local result is clear from
Theorem~\vref{theorem:ReductiveGstructuresHartogsExtension}.
Take a $G$-structure
$ \phi : E \to FT\left(M \backslash S\right)$
on $M \backslash S$, with~$M$ an $n$-fold,
and $S$ a subset of \cct.
In local coordinates $z^1, \dots, z^n$ on~$M$ near a point $s \in S$,
points of $FTM$ look like $\left(z,u\right)$
with $z \in \C{n}$ and $u \in \GL{n,\C{}}$.
The map $\Gromov{\phi}$ is $\Gromov{\phi}(z,u)=uf(z)$ for some
map $f : M \backslash S \to \GL{n,\C{}}/G$. Extend $f$ to a holomorphic
map $f : M \to \GL{n,\C{}}/G$. Now let $E'$ be the set of points
of the form $(z,u) \in FTM$ so that $uf(z)=G \in \GL{n,\C{}}/G$.
This bundle $E' \to M$ extends $E$. By
Theorem~\vref{theorem:GeometricStructure},
the $G$-structure extends holomorphically.

Next suppose that $\GL{n,\C{}}/G$ is not a
\Thullen extension
target for $n$-folds. Pick
an $n$-dimensional complex manifold $M$ and
a holomorphic map
$f : M \backslash S \to \GL{n,\C{}}/G$
which does not extend to $M$.
There must be some point $s \in S$
near which $f$ doesn't extend holomorphically.
We can replace
$M$ by any neighborhood of $s$, so we can
assume that $M$ is a ball $B$ in $\C{n}$.
Think of $\GL{n,\C{}} \to \GL{n,\C{}}/G$ as
a holomorphic principal right $G$-bundle.
Let $E \to B \backslash S$ be the pullback bundle $E = f^* \GL{n,\C{}}$.
By def\/inition, $E$ is a principal right $G$-subbundle of
$\left(B \backslash S\right) \times
\GL{n,\C{}}=FT\left(B \backslash S\right)$,
hence a $G$-structure on $B \backslash S$.
Suppose that this $G$-structure~$E$ extends
holomorphically to a $G$-structure~$E'$ on $B$.
If need be, we replace $B$ by a smaller
ball around $s$ on which the bundle
$E'$ is trivial, $E'=B \times G$.
Then we can take the section $B \times \left\{1\right\}$
of $E'$, and map it to $\GL{n,\C{}}$ and
then quotient by $G$ to extend $f$ to
a holomorphic map $f : B \to \GL{n,\C{}}/G$.
\end{proof}

\begin{Example}[Hwang and Mok \cite{Hwang/Mok:1998}]\label{example:ReductiveAlgebraicG}
By Lemma~\vref{lemma:ClassicalInvariantTheory},
if $G$ is a reductive algebraic group, then $\GL{n,\C{}}/G$ is an af\/f\/ine
variety, so holomorphic $G$-structures extend
across subsets of \cct
and extend from
domains in Stein manifolds
to their envelopes of holomorphy.
\end{Example}

\begin{Example}
Again consider holomorphic Riemannian metrics: $G=\SO{n,\C{}}$.
Then
\[
\GL{n,\C{}}/\SO{n,\C{}}=X \backslash (f=0),
\]
with $X$ the set of complex quadratic forms, and $f$ the determinant.
Therefore $\SO{n,\C{}}$-structures (holomorphic Riemannian
metrics) extend across subsets of \cct,
and extend from domains in Stein manifolds to their envelopes of holomorphy.
\end{Example}

\begin{Example}
Consider almost symplectic structures, $G=\Symp{2n,\C{}}$.
Because $G$ is reductive, almost symplectic
structures extend across subsets of \cct. Alternately, pick
a complex volume form $\Omega \in \Lm{2n}{\C{2n}}^*$, and def\/ine
$f : \Lm{2}{\C{2n}}^* \to \C{}$ by $f(\alpha)=\alpha^n/\Omega$.
Then $\GL{2n,\C{}}/\Symp{2n,\C{}} \backslash (f=0) \subset \Lm{2}{\C{2n}}^*$
is the set of symplectic forms.
Therefore $\Symp{2n,\C{}}$-structures (holomorphic almost symplectic
structures) extend across subsets of \cct,
and extend from domains in Stein manifolds to their envelopes of holomorphy.
\end{Example}

\begin{Example}
Given a 3-form $\sigma$ on $\C{7}$,
def\/ine $B_{\sigma} : \C{7} \otimes \C{7} \to \Lm{7}{\C{7}}^*$
by $B_{\sigma}(u,v)=\frac{1}{6} \left(u \hook \sigma\right)
\wedge \left( v \hook \sigma\right) \wedge \sigma$.
Pick some nonzero $\Omega \in \Lm{7}{\C{7}}^*$.
Let $f(\sigma)=\det \left( B_{\sigma}/\Omega \right)$.
Say that $\sigma$ is \emph{nondegenerate} if $f(\sigma) \ne 0$, i.e.
$B_{\sigma}/\Omega$ is a nondegenerate
quadratic form. For instance, if we write $dz^{ij}$ to
mean $dz^i \wedge dz^j$, etc., then
the 3-form
\[
\sigma_0 =
dz^{123} + dz^{145} + dz^{167} + dz^{246} - dz^{257} - dz^{347} - dz^{356}
\]
is nondegenerate. The degenerate forms clearly form an af\/f\/ine analytic
hypersurface $(f \ne 0)$ inside the space of $3$-forms.

It turns out (see Bryant~\cite{Bryant:1987}) that the nondegenerate
$3$-forms are precisely the orbit of $\sigma_0$ under
$\GL{7,\C{}}$-action in $\Lm{3}{\C{7}}^*$.
Moreover the stabilizer of $\sigma_0$ is the exceptional simple
Lie group $G_2$. Therefore $\GL{7,\C{}}/G_2=\Lm{3}{\C{7}}^* \backslash (f=0)$
is a \Thullen and Hartogs extension target. So holomorphic $G_2$-structures
(even with torsion) extend across subsets of \cct
and extend from
domains in Stein manifolds to their envelopes of holomorphy.
\end{Example}

\begin{Example}\label{example:web}
A \emph{web} on a surface $M$ is a choice
of three nowhere tangent foliations by curves.
We will see that webs extend
across punctures.
At each point $m \in M$, there is a linear
isomorphism $T_m M \to \C{2}$ taking the tangent
lines of the curves to the horizontal axis,
vertical axis and diagonal. This linear isomorphism
is unique up to rescaling. If we let $E$ be the
set of all such linear isomorphisms at all points
of $M$, then $E \to M$ is a $\C{\times}$-structure.
A web is therefore a $\C{\times}$-structure.
Clearly $\GL{2,\C{}}/\C{\times}=\PSL{2,\C{}}$ is covered
by $\SL{2,\C{}}$, an af\/f\/ine variety so a \Thullen extension
target. Therefore holomorphic webs extend
across punctures, and extend from
domains in Stein manifolds to their
envelopes of holomorphy. This is surprising because foliations
by curves in a surface need not extend across punctures.

It turns out that a web is equivalent to a
$\Proj{1} \times \Proj{1} \backslash \text{diagonal}$-geometry;
see Example~\vref{example:RiemannCartanSurfaces}. Webs can also
be described as a certain type of $\C{2}$-geometry, i.e.\
as modelled on $G/H=\C{2}$ for certain groups $H$ and $G$,
which we leave to the reader to work out.
\end{Example}

\begin{Example}
There is a cubic form on $\C{27}$ whose stabilizer
is $E_6$. Since $E_6$ is a reductive algebraic group,
holomorphic $E_6$-structures
(i.e. symmetric cubic forms in tangent spaces, with
stabilizer isomorphic to $E_6$) extend across subsets of \cct
and extend from
domains in Stein manifolds to their envelopes of holomorphy.
\end{Example}

\begin{Example}
The space of solutions of certain scalar complex analytic ordinary dif\/fe\-ren\-tial
equations bears a canonical $\GL{2,\C{}}$-structure,
for a particular embedding $\GL{2,\C{}} \subset \GL{5,\C{}}$;
see Doubrov \cite{Doubrov:2007}, Dunajski and Tod \cite{Dunajski/Tod:2005}
or Godli\'nski and Nurowski \cite{Godlinski/Nurowski:2007}.
Because $\GL{2,\C{}}$ is a reductive algebraic group,
these structures extend over subsets of \cct
and extend from
domains in Stein manifolds to their envelopes of holomorphy.
\end{Example}

\begin{Example}\label{example:AlmostProduct}
 An \emph{almost product structure} is a $G$-structure
where
\[
G = \GL{k,\C{}} \times \GL{n-k,\C{}} \subset \GL{n,\C{}}.
\]
Equivalently, an almost product structure is a pair of
complementary transverse plane f\/ields.
Clearly $G$ is a reductive algebraic group, so almost product structures
extend across subsets of \cct
and extend from
domains in Stein manifolds to their envelopes of holomorphy.
\end{Example}

\begin{Example}[Ivashkovich \protect{\cite[Proposition~3, p.~196]{Ivashkovich:1985}}]
If $X$ and $Y$ are \Thullen extension targets for local biholomorphisms,
then so is $X \times Y$.
\end{Example}

In fact, a much stronger result is true, which
does not apparently follow from our theorems above.
\begin{Theorem}[Ivashkovich \protect{\cite[Proposition~4, p.~196]{Ivashkovich:1985}}]
If the base and fiber of a holomorphic fibration
are Hartogs extension targets, then so is the total space.
\end{Theorem}

\begin{Example}
Consider a Hopf f\/ibration: let
$M = \left(\C{3} \setminus 0\right) / \left(z \sim 2z\right)$
and map $F : M \to \Proj{2}$, $F(z)= \C{}z \in \Proj{2}$.
The map $F$ is a f\/ibration $E \to M \to \Proj{2}$ with f\/iber
$E=\left(\C{} \setminus 0\right) / \left(z \sim 2z\right)$
an elliptic curve.
By Theorem~\vref{theorem:KahlerHomogeneous}, the base is a
Hartogs and \Thullen target for local biholomorphisms.
By Lemma~\ref{lemma:coverings}, the f\/iber
is a Hartogs and \Thullen target. The total space
admits the obvious
map $\C{3} \backslash 0 \to M$ which does not extend
to $\C{3}$. Therefore the total space is neither
a Hartogs nor a \Thullen target, even for local
biholomorphisms, . Similarly we can construct the obvious bundle
$F \times F : M \times M \to \Proj{2} \times \Proj{2}$
with f\/iber $E \times E$, and
map $\left(\C{3} \backslash 0\right)
\times \left(\C{3} \backslash 0\right) \to M \times M$.
The map doesn't extend to $\C{6}$, i.e.\
over the complex codimension~3 subset $0 \times \C{3} \cup \C{3} \times 0$.
The reader can build many more examples along the same lines
in various dimensions. It might be
signif\/icant that $M$ is not K\"ahler.
\end{Example}

\begin{Conjecture}\label{conjecture:LocalBiholProduct}
Consider a holomorphic fibration whose total
space is K\"ahler. Suppose that the base and fiber
are Hartogs extension targets for local biholomorphisms.
Then so is the total space.
\end{Conjecture}

\subsection{Relative extension problems for f\/irst order structures}

\begin{Definition}
Suppose that $G_0$ and $G_1$ are complex Lie groups,
and that $\rho_0 : G_0 \to G_1$ and $\rho_1 : G_1 \to \GL{n,\C{}}$
are Lie group morphisms. We then treat $\C{n}$ as both
a~$G_0$-module and a~$G_1$-module, using representations
$\rho_1 \rho_0$ and $\rho_1$ respectively.
Every $G_0$-structure $\phi : E_0 \to FTM$ induces a $G_1$-bundle
$E \times_{G_0} G_1 \to M$.
We def\/ine $\phi : E \times G_1 \to FTM$ by
$\phi\left(e,g_1\right)=\rho_1(g_1)^{-1} \phi(e)$. This map
clearly descends to $E \times_{G_0} G_1$, giving
the \emph{induced $G_1$-structure}.
\end{Definition}

\begin{Theorem}\label{theorem:Relative}
 Suppose that $G_0$ and $G_1$ are complex Lie groups,
and that $\rho_0 : G_0 \to G_1$ and $\rho_1 : G_1 \to \GL{n,\C{}}$
are Lie group morphisms, with closed images.
If a $G_0$-structure extends across a subset
then the induced $G_1$-structure does as well.
Suppose that $\rho_0$ has discrete kernel
and that holomorphic maps to $G_1/\rho_0\left(G_0\right)$ extend
across subsets of \cct.
Then any $G_0$-structure extends over a subset of \cct
just when its induced $G_1$-structure extends over the same subset.
\end{Theorem}
\begin{proof}
Combine Theorem~\vref{theorem:GeometricStructure} and Lemma~\vref{lemma:Relative}.
\end{proof}

\begin{Example}
Holomorphic spin structures, $G=\Spin{n,\C{}}$,
extend across subsets of \cct, because the
induced $\SO{n,\C{}}$-structures do.
\end{Example}

\subsection{Torsion and invariant relations}\label{section:FullTorsion}

\begin{Definition}
Suppose that $G$ is a complex Lie group with Lie algebra $\mathfrak{g}$
and $V_0$ is a complex $G$-module, say with representation
$\rho : \mathfrak{g} \to \gl{V_0}$.
Let $\delta : \gl{V_0} \otimes V_0^* \to V_0 \otimes \Lm{2}{V_0}^*$
by $\delta\left(A \otimes \xi \right)(v,w)=\rho(A)(v)w-\rho(A)(w)v$.
We def\/ine $\prol{\mathfrak{g}}$ and $\Cohom{0,2}{\mathfrak{g}}$ to be the
kernel and cokernel of $\delta$:
\[
\xymatrix{ 0 \ar[r] & \prol{\mathfrak{g}} \ar[r] & \mathfrak{g} \otimes V^*
  \ar[r]^{\delta} & V_0 \otimes \Lm{2}{V_0}^* \ar[r]^{[]}
& \Cohom{0,2}{\mathfrak{g}} \ar[r] &
  0  }
\]
denoting the quotient map to the cokernel as $t \mapsto [t]$.

Suppose that $E {\to} FTM$ is a $G$-structure.
Denote the projection $FTM {\to} M$ as  $\pi : FTM {\to} M$ .
We def\/ine the \emph{soldering $1$-form}
$\sigma \in \nForms{1}{FTM}$ by $v \hook \sigma_u = u(\pi'(u)v)$.
We will also denote the pullback of $\sigma$ to $E$ as $\sigma$.

For each $A \in \mathfrak{g}$, we write
the associated inf\/initesimal generator of the right
$G$-action on~$E$ as $\vec{A}$.
Denote the projection $E \to M$ as $\pi_E : E \to M$.
For each open set $U \subset M$,
a \emph{pseudoconnection} for $E$ over $U$
is a choice of 1-form
$\gamma \in \nForms{1}{\pi_E^{-1}U} \otimes \mathfrak{g}$ so that
$\vec{A} \hook \gamma = A$ for all $A \in \mathfrak{g}$.
The local existence of a pseudoconnection is obvious,
because $E$ is locally trivial.
Any two pseudoconnections $\gamma$ and $\gamma'$
def\/ined over the same open set $U \subset M$ dif\/fer
by $\gamma'-\gamma = Q \sigma$, where $Q : \pi_E^{-1} U \to \mathfrak{g} \otimes V_0^*$.
Moreover, any choice of such a function $Q$ yields
a new pseudoconnection $\gamma'=\gamma+Q \sigma$.

The \emph{torsion} of a pseudoconnection $\gamma$ is the function
$t : \pi_E^{-1} U \to V_0 \otimes \Lm{2}{V_0}^*$ for which
\[
 d \sigma + \gamma \wedge \sigma = \frac{1}{2} t \sigma \wedge \sigma.
\]
If we change the pseudoconnection to $\gamma'=\gamma+Q \omega$, then
we change the torsion to $t'=t+\delta Q$. Therefore
the \emph{intrinsic torsion} $[t] : E \to \Cohom{0,2}{\mathfrak{g}}$
of the $G$-structure is well def\/ined globally.
A~\emph{torsion relation} is a $G$-invariant subset of
$\Cohom{0,2}{\mathfrak{g}}$. An \emph{integral
structure} of a torsion relation is a f\/irst
order structure whose torsion lies in the torsion relation.

For simplicity, assume that $\rho : G \to \GL{V_0}$ has discrete kernel.
Consider the map
$[\delta] : \left(\gl{V_0}/\mathfrak{g}\right) \otimes V_0^* \to \Cohom{0,2}{\mathfrak{g}}$,
given on $S \in \gl{V_0} \otimes V_0^*$ by $[\delta](S)=[\delta S]$.
Given a torsion relation $R \subset \Cohom{0,2}{\mathfrak{g}}$,
let $R'$ be the preimage of $R$ under $[\delta]$, and
call $R'$ the \emph{induced invariant relation}
for maps to $\GL{V_0}/G$.
\end{Definition}

\begin{Theorem}\label{theorem:InvariantRelationsAndTorsionRelations}
 Suppose that $G \to \GL{n,\C{}}$ is a representation with discrete kernel
and closed image.
Let $R$ be a torsion relation for $G$, and $R'$ the induced invariant
relation for maps to $\GL{n,\C{}}/G$.
Integral maps of $R'$ extend across
subsets of \cct just when integral structures of $R$
extend across subsets of \cct.
\end{Theorem}
\begin{proof}
Suppose that $M$ is a complex manifold, $S \subset M$ a subset of \cct,
and $E \to FT(M \backslash S)$ is a $G$-structure which is an
integral structure for $R$.
We only need to extend $E$ locally, so we can assume that $M$ is an
open subset of $\C{n}$, and so $FTM = M \times \GL{n,\C{}}$. By Theorem~\vref{theorem:Relative},
we can assume that $G \subset \GL{n,\C{}}$. Therefore $E$ is
a subbundle of the trivial bundle $\left(M \backslash S\right) \times \GL{n,\C{}}$,
and is therefore determined as the pullback bundle $f^*\GL{n,\C{}}$
of the bundle $\GL{n,\C{}} \to \GL{n,\C{}}/G$ by a map
$f : M \backslash S \to \GL{n,\C{}}/G$.
So if we can extend $f$ across $S$, then we can extend $E$ across $S$ to a~principal right $G$-bundle $E \to M$, and therefore by
Theorem~\vref{theorem:GeometricStructure} the $G$-structure
extends across $S$. Therefore we need only prove that
$f$ extends across $S$.

Take any local section $u$ of $E \to M \backslash S$. This section is
then a matrix $u : \text{open} \subset M \backslash S \to \GL{n,\C{}}$. Then clearly $f=uG$.
We need to compute $f' : T_m M \to T_{f(m)} \left(\GL{n,\C{}}/G\right)$.
If we let $\pi : \GL{n,\C{}} \to \GL{n,\C{}}/G$ be the bundle
map, then $f=uG=\pi u$, so $f'(m)=\pi'(u(m))u'(m)$.
Denote left translation by any matrix $g \in \GL{n,\C{}}$ as $L_g$.
Let us see why $f$ is an integral map of $R'$.
We need to compute $f^{(1)} : FTM \times_M f^*\GL{n,\C{}} \to
\left(\mathfrak{\gl{n,\C{}}}/\mathfrak{g}\right) \otimes \C{n*}$.
But $FTM=M \times \GL{n,\C{}}$, so $FTM \times_M f^*G=M \times \GL{n,\C{}} \times G$,
identif\/ied by taking $(m,h,g) \in M \times \GL{n,\C{}} \times G
\to (m,h,u(m)g) FTM \times_M f^*G$. So
$f^{(1)}(m,h,g)=\left(L_{u(m)g}^{-1}\right)'(f(m)) f'(m)h^{-1}$.
In order to test if $f$ is an integral map of $R'$,
it suf\/f\/ices to take $g$ and $h$ to be the identity:
\begin{gather*}
\left(L_u^{-1}\right)'f'(m)
=
\left(L_u^{-1}\right)'\pi'(u(m))u'(m)
=
\left(L_u^{-1}\pi\right)'(u(m))u'(m)
 =
\left(\pi L_u^{-1}\right)'(u(m))u'(m)
\\
\phantom{\left(L_u^{-1}\right)'f'(m)}{}  =
\pi'(m) \left(L_u^{-1}\right)'(u(m))u'(m)
 =
\pi'(m) u^{-1} \, du
 =
u^{-1} \, du + \mathfrak{g}.
\end{gather*}
So it suf\/f\/ices to show that $u^{-1} \, du + \mathfrak{g} \in
\nForms{1}{E} \otimes \left(\gl{n,\C{}}/\mathfrak{g}\right)$
is a 1-form valued in $R'$.

We can work entirely locally on $M$, so we can assume that $M$ is a domain in $\C{n}$
with coordinates $z^1, z^2, \dots, z^n$. Then $\sigma=u^{-1} \, dz$, and
so $d \sigma = - \left( u^{-1} \, du \right) \wedge \sigma$.
But we also have $d \sigma = - \gamma \wedge \sigma + \frac{1}{2} t \sigma \wedge \sigma$.
By Cartan's lemma, we can write $u^{-1} \, du = \gamma + \frac{1}{2} t \sigma + \frac{1}{2} Q \sigma$
with $Q \in \C{n} \otimes \Sym{2}{\C{n}}^*$. So
\[
 u^{-1} \, du + \mathfrak{g} = \frac{1}{2} t \sigma + \frac{1}{2} Q \sigma + \mathfrak{g}.
\]
Clearly $\delta(t+Q)=t$, so $[\delta(t+Q)]=[t]$. Because $[t] \in R$,
we must have (independent of choices made of trivializations and
pseudoconnections) $t+Q + \mathfrak{g} \otimes V_0^* \in R'$,
so that $u^{-1} \, du + \mathfrak{g}$ is valued in $R'$.
Therefore if integral maps of $R'$ extend
across subsets of \cct, then integral structures
of $R$ extend as well.

Suppose that $f : M \backslash S \to \GL{n,\C{}}/G$ is an
integral map of $R'$. Take the bundle $\GL{n,\C{}} \to \GL{n,\C{}}/G$
and pullback to a $G$-bundle $E=f^* \GL{n,\C{}}$.
In order to extend $f$ across $S$, it suf\/f\/ices to do so
locally, so we can assume that $M$ is an open subset
of $\C{n}$, and so $FTM$ is holomorphically
trivial, $FTM = M \times \GL{n,\C{}}$. We can then
map $E \to FTM$ by the identity map, a $G$-structure. The torsion
is then clearly in $R$ by the same arguments as above.
\end{proof}

The reader can see how to prove a relative version of the previous theorem.

\begin{Conjecture}\label{conjecture:InvariantRelationsAndTorsionRelations}
 Suppose that $G \to \GL{n,\C{}}$ is a representation with discrete kernel
and closed image.
Let $R$ be a torsion relation for $G$, and $R'$ the induced invariant
relation for maps to $\GL{n,\C{}}/G$.
Integral maps of $R'$ extend from any domain $M$ in any Stein manifold
to the envelope of holomorphy of $M$
just when integral structures of $R$ extend from any domain $M$ in any Stein manifold
to the envelope of holomorphy of $M$.
\end{Conjecture}

\begin{Example}\label{example:ContactRelation}
 We have seen that holomorphic contact structures extend across
subsets of \cct. We can consider a contact structure
as a hyperplane f\/ield, i.e.\ a $G$-structure where $G$ is the
subgroup of $\GL{2n+1,\C{}}$ preserving the hyperplane
$\C{2n} \subset \C{2n+1}$.
The torsion module is
$\Cohom{0,2}{\mathfrak{g}}
= \Lm{2}{\C{2n}}^* \otimes \left(\C{2n+1}/\C{2n}\right)$.
If $t$ is the torsion of a $G$-structure, then
$t^n \in  \Lm{2n}{\C{2n}}^* \otimes \left(\C{2n+1}/\C{2n}\right)^{\otimes n}$.
A $G$-structure is a contact structure just when $t^n \ne 0$.
Let $R \subset \Cohom{0,2}{\mathfrak{g}}$
be the open set of $t \in \Cohom{0,2}{\mathfrak{g}}$ so that $t^n \ne 0$.
Looking at the relevant matrices, it is easy to see that
$\gl{2n+1,\C{}}/\mathfrak{g}=\C{2n*} \otimes \left(\C{2n+1}/\C{2n}\right)$.
The induced invariant relation $R' \subset \left(\gl{2n+1,\C{}}/\mathfrak{g}\right) \otimes \C{2n+1*}$
is the set of tensors $L \in \C{2n*} \otimes \left(\C{2n+1}/\C{2n}\right) \otimes \C{2n+1*}$
for which the induced tensor $L \in \C{2n*} \otimes \left(\C{2n+1}/\C{2n}\right) \otimes \C{2n*}$
has symplectic antisymmetrization in the obvious indices.

We conclude that if $M$ is any complex manifold of
dimension $2n+1$, and $S \subset M$ any subset of \cct,
then any holomorphic map $f : M \backslash S \to \Proj{2n}$
satisfying this invariant relation $R'$
extends across $S$. We can compute the resulting
dif\/ferential relation on $f$ easily. In local coordinates, write
\[
 f(z)=
\begin{bmatrix}
 1 \\
f_1(z) \\
\vdots \\
f_{2n}(z)
\end{bmatrix}.
\]
The $f$ is an integral map for $R'$ just when
\[
\Omega = \left( \pd{f_j}{z_k} - \pd{f_k}{z_j}\right) dz^j \wedge dz^k
\]
satisf\/ies $\Omega^n \ne 0$.
By Theorem~\vref{theorem:InvariantRelationsAndTorsionRelations},
the relation $R'$ is a \Thullen extension relation,
because contact structures extend holomorphically
across subsets of \cct.
Indeed just looking at this relation, it is clearly a \Thullen extension
relation. It is not so obvious that $R'$ is invariant under projective
transformations.
This relation $R'$ is the only known \Thullen extension relation
for maps $M^{2n+1} \to \Proj{2n}$.
\end{Example}

\subsection{Harmless reductions of f\/irst order structures}\label{section:HarmlessReductions}

A simple trick allows us to ef\/fectively reduce
f\/irst order structures for the sake of
solving the \Thullen extension problem.
Suppose that $G \subset L \subset \GL{n,\C{}}$,
each a closed subgroup of the next.
Suppose that $L_0 \subset L$ is another closed
subgroup and that $G$ acts transitively on
$L/L_0$. Then we will say that
the pair of subgroups $L_0 \subset L$ are
\emph{harmless} for $G$. We will then
let $G_0 \subset G$ be the subgroup
f\/ixing the identity coset of $L/L_0$, so $G/G_0=L/L_0$,
and $G_0 = G \cap L_0$.

For example, $L=\GL{n,\C{}}$ and $L_0=\SL{n,\C{}}$
are harmless for any $G \subset \GL{n,\C{}}$
whose identity component is not contained in $\SL{n,\C{}}$.

We can reduce the \Thullen extension
problem for $G$-structures to the \Thullen
extension problem for $G_0$-structures.
\begin{Lemma}
Suppose that $M$ is a complex manifold
and that $S \subset M$ is a
subset of \cct. Suppose that $E \to FT\left(M \backslash S\right)$
is a $G$-structure.
Take any harmless pair $L_0 \subset L$ for $G$.
Then $E \to FT\left(M\backslash S\right)$
induces an $L$-structure
$E \times_G L$ as above.
Suppose that the $L$-structure extends
across $S$. Then $M$ is covered by open sets
$U_a$ so that $\left.E\right|_{U_a \backslash S}$
admits a holomorphic reduction
to a $G_0$-structure $E_a$. The
$G$-structure $E$ extends across $S$
just when every one of these $G_0$-structures $E_a$
extends across $U_a \cap S$.
\end{Lemma}
\begin{proof}
Suppose that $E \times_G L \to M \backslash S$ extends
to some $L$-structure $E' \to M$.
Pick any open set $U \subset M$ on which $E'$ admits
an $L_0$-reduction, say $E_0' \to U$. Then let $E_0=E \cap E_0'$.

Let us see why $E_0 \to U \backslash S$ is a $G_0$-structure. Take any open
set $U_0 \subset U$ on which $E$ is holomorphically trivial,
say $\left.E\right|_{U_0} = U_0 \times G$. Then $E'= U_0 \times L$.
If we shrink $U_0$ we can arrange that the $L_0$-structure is also
trivial, $E_0' = U_0 \times L_0$, mapped to $E'=U_0 \times L$
by some bundle map $\phi\left(m,\ell_0\right)=r_{\ell_0} \phi(m)$,
for any $m \in U_0$ and $\ell_0 \in L_0$, for
some holomorphic map $\phi : U_0 \to L$. So then $E_0$ is the set of pairs
$\left(m,\ell_0\right)$ for which $r_{\ell_0} \phi(m) \in G$.
Clearly $E_0$ is acted on freely by $G_0$, since the
equation $r_{\ell_0} \phi(m) \in G$ is $G_0$-equivariant. We only
need to show that $G_0$ acts transitively on the f\/ibers
of $E_0 \to U_0$. Consider two points lying in the same
f\/iber on $E_0 \to U_0$, say $\left(m,\ell_0\right)$ and
$\left(m,\ell_1\right)$. So $r_{\ell_0} \phi(m) \in G$
and $r_{\ell_1} \phi(m) \in G$. Let $\ell = \phi(m)$.
We have $\ell \ell_0 \in G$ and $\ell \ell_1 \in G$.
So $\left(\ell \ell_0 \right)^{-1} \ell \ell_1 \in G$.
So $\ell_0^{-1} \ell_1 \in G$ and therefore $\ell_0^{-1} \ell_1 \in G_0$.
We see that $E_0 \to U_0$ is a holomorphic $G_0$-subbundle
of $E$.
\end{proof}

\begin{Example}
 If $G \subset \GL{n,\C{}}$ has identity
component not contained in $\SL{n,\C{}}$,
then any $G$-structure can be locally
reduced to a $G_0$-structure,
$G_0=G \cap \SL{n,\C{}}$, and the $G$-structure
extends across a subset of \cct just
when all of the $G_0$-structures do;
roughly speaking we can pick a local
holomorphic volume form.
\end{Example}

Once we have achieved a harmless reduction,
we can apply Cartan's method of
equivalence (see Gardner \cite{Gardner:1989},
Ivey and Landsberg \cite{Ivey/Landsberg:2003}
or Sternberg \cite{Sternberg:1983})
to this $G_0$-structure, to try
to reduce it further, if possible.

\begin{Example}\label{example:scalarConservationLaw}
 For this example (but not for any subsequent
part of this paper), we will expect the reader
to be conversant with Cartan's method
of equivalence. Bryant, Grif\/f\/iths and
Hsu \cite{Bryant/Griffiths/Hsu:1995}
constructed out of any scalar conservation law
an equivalent f\/irst order structure. We will see
why holomorphic scalar conservation laws
extend across subsets of \cct.

Their $G$-structure
has structure equations (in a slight alteration
of their notation)
\[
 d
\begin{pmatrix}
 \omega^1 \\
\omega^2 \\
\omega^3
\end{pmatrix}
= -
\begin{pmatrix}
 2 \phi & 0 & 0 \\
0 & \phi & 0 \\
0 & \mu & - \phi
\end{pmatrix}
\wedge
\begin{pmatrix}
 \omega^1 \\
\omega^2 \\
\omega^3
\end{pmatrix}
+
\begin{pmatrix}
K \omega^2 \wedge \omega^3 \\
\omega^1 \wedge \omega^3  \\
0
\end{pmatrix}.
\]
The structure group $G$ is the group of
matrices of the form
\[
 \begin{pmatrix}
  g^2 & 0 & 0 \\
  0   & g & 0 \\
  0   & h & g^{-1}
 \end{pmatrix}
\]
for any nonzero real number $g$ and any real number $h$.

Complexify: consider holomorphic scalar conservation
laws, so our group $G$ has $g$ and $h$ complex.
Picking a local holomorphic volume form,
as described in the previous example,
we can harmlessly reduce to the group of
matrices of the form
\[
 \begin{pmatrix}
  1 & 0 & 0 \\
  0   & \pm 1 & 0 \\
  0   & h & \pm 1
 \end{pmatrix}.
\]
This structure group is still not a reductive
algebraic group, but we can apply Cartan's
method of equivalence to the reduced structure.
The structure equations are identical, but
$\phi$ becomes semibasic, i.e.
\[
 \phi = -t_1 \, \omega^1 - t_2 \, \omega^2 - t_3 \, \omega^3
\]
for some holomorphic functions $t_1$, $t_2$, $t_3$ on
the total space of each $G_0$-structure. (The minus signs
are for convenience in the following calculations.)

Take exterior derivatives of all of the structure equations,
to see that $dt_1 = \mu$ on the f\/ibers of the
total space of each $G_0$-structure. So
$t$ translates under $G_0$-action, and therefore
the set $B_1=\left(t_1=0\right)$ is a $G_1$-subbundle, where
$G_1$ is the group of matrices of the form
\[
 \begin{pmatrix}
  1 & 0 & 0 \\
  0   & \pm 1 & 0 \\
  0   & 0 & \pm 1
 \end{pmatrix}.
\]
The structure group is now a reductive algebraic group,
so we can extend each $G_1$-structure across the
subset of \cct, and therefore extend the
original $G$-structure.
\end{Example}

\begin{Example}
For Engel 2-plane f\/ields on 4-folds,
harmless reduction doesn't
help solve the \Thullen extension
problem.
The structure group of an Engel 2-plane f\/ield
(see Ehlers et al.~\cite{Ehlers/Koiller/Montgomery/Rios:2005})
can only extend to a reductive algebraic
group $L$ by taking $L=\GL{4,\C{}}$.
Then the only choice of $L_0$ is
$\SL{4,\C{}}$. This harmless reduction reduces the Engel
structure by picking a volume form. But the
symmetry pseudogroup of the standard Engel
2-plane f\/ield on the standard 2-jet
bundle together with the standard volume form is
given by prolongations of maps $X=X(x), Y=Y(x,y)=Y\left(x,y_0\right)+y \, X'(x)^{1/3}$.
This pseudogroup is still inf\/inite dimensional,
so the structure is still of inf\/inite type,
and the structure group is not a reductive algebraic group.
\end{Example}

\begin{Example}
The author hopes that harmless reduction
might solve the \Thullen extension problem
for Clelland's $G$-structure
associated to a parabolic partial dif\/ferential
equation for one function of $1+2$ variables
to prove that these equations extend
across subsets of \cct;
see Clelland \cite{Clelland:1997}. Note that such
a result would only extend the dif\/ferential equation, and
not its solutions.
\end{Example}

\begin{Example}
 If $G \subset \GL{n,\C{}}$ preserves
a hyperplane $P$, then a holomorphic $G$-structure induces
a holomorphic hyperplane f\/ield. If that hyperplane f\/ield
is a holomorphic contact structure, then
that contact structure extends across subsets of \cct.
If $G$ does not preserve a 1-form
vanishing on the hyperplane $P$, then
picking a local choice of contact form
is a harmless reduction; this was our method
in Theorem~\vref{theorem:ContactStructuresExtend}.
\end{Example}

We will also apply harmless reduction to holomorphic
parabolic geometries in Theorem~\vref{theorem:ParabolicGeometriesExtend}.
There is an obvious analogue of harmless reduction
for maps to homogeneous spaces, following
the ideas of Theorem~\vref{theorem:InvariantRelationsAndTorsionRelations},
which we leave the reader to explore.

\begin{Remark}
We leave it to the reader to generalize
harmless reduction to the Hartogs
extension problem for a domain $M$
in a Stein manifold, as long as
$M$ bears a holomorphic volume form
(in particular for Riemann domains).
\end{Remark}

\section{Higher order structures}\label{section:HigherOrderStructures}

\begin{Definition}
Fix a complex manifold $M$, a vector space $V_0$ with $\dim V_0=\dim M$,
and take $FTM$ the $V_0$-valued frame bundle.
Let $\pi : FTM \to M$ be the bundle mapping.
We def\/ine a~1-form $\sigma$ on $FTM$, called the \emph{soldering form},
by $v \hook \sigma_{u} = u\left(\pi'(u)v\right)$.
Suppose that $E \to FTM$ is a $G$-structure.
Let $\mathfrak{g}$ be the Lie algebra of $G$.
For any element $A \in \mathfrak{g}$, write
$\vec{A}$ for the vector f\/ield on $E$ generating
the right $H$-action, i.e.
\[
\left.\frac{d}{dt} r_{e^{tA}} \right|_{t=0} = \vec{A}.
\]
We will always also denote the pullback of $\sigma$
on $E$ as $\sigma$. A \emph{pseudoconnection $1$-form}
at a point $e \in E$ is a 1-form $U \in T^*_e E \otimes \mathfrak{g}$
so that $\vec{A} \hook U=A$.
\end{Definition}
Pseudoconnection 1-forms exist at each point $e$ of $E$.
The set of all pseudoconnection
1-forms at all points of the total space $E$ of a $G$-structure
form a principal right $\mathfrak{g} \otimes V_0^*$-bundle over $E$,
under the action $r_{A \otimes \xi}U=U-(A \sigma)\wedge (\xi \sigma)$.
The reader can consult Gardner \cite{Gardner:1989}
or Ivey and Landsberg \cite{Ivey/Landsberg:2003}.

\begin{Definition}
A \emph{torsion function} on a $G$-structure $E \to FTM$ is a
holomorphic function $t : E \to \Lm{2}{V_0}^* \otimes V_0$ so that,
\begin{enumerate}\itemsep=0pt
\item $t$ is $G$-equivariant,
 \item at each point $e \in E$, there is a pseudoconnection
1-form $U$ so that
$d \sigma + U \wedge \sigma = \frac{1}{2} t \sigma \wedge \sigma$.
\end{enumerate}
\end{Definition}
Not every $G$-structure admits a torsion function in this sense,
but the most important examples of $G$-structures do;
see Gardner \cite{Gardner:1989}
or Ivey and Landsberg \cite{Ivey/Landsberg:2003}.

\begin{Definition}
 If $\rho : G \to \GL{V_0}$ is a holomorphic representation of a complex
Lie group $G$ with Lie algebra $\mathfrak{g}$, let
$\delta : \mathfrak{g} \otimes V_0^* \to V_0 \otimes \Lm{2}{V_0}^*$
be def\/ined by $\delta\left(A \otimes \xi\right)(v_1,v_2)=
\left(Av_1\right)\xi\left(v_2\right)-\left(Av_2\right)\xi\left(v_1\right)$.
We def\/ine $\prol{\mathfrak{g}}=\ker \delta$.
\end{Definition}

\begin{Definition}\label{definition:Prolongation}
If $E \to FTM$ is a $G$-structure with torsion function $t$,
then def\/ine the \emph{prolongation} $\prol{E}$ of $E$ (with
respect to $t$) to be the bundle of all pairs
$\left(e,U\right)$ of points $e \in E$ and
pseudoconnection 1-forms $U$ at $e$
with $d \sigma + U \wedge \sigma = \frac{1}{2} t \sigma \wedge \sigma$.
\end{Definition}
It is easy to check that $\prol{E} \to E$ is a principal
right $\prol{\mathfrak{g}}$-bundle (a subbundle of the bundle
of pseudoconnection 1-forms), with $\prol{\mathfrak{g}}$ acting
as a subgroup of $\mathfrak{g} \otimes V_0^*$.

\begin{Example}
Suppose that $\pi : E \to M$ is a $G/H$-geometry with Cartan connection
$\omega$. Let $\sigma = \omega+\mathfrak{h} \in \nForms{1}{E} \otimes
\left(\mathfrak{g}/\mathfrak{h}\right)$.
There is a unique function
$K : E \to \Lm{2}{\mathfrak{g}/\mathfrak{h}}^* \otimes \mathfrak{g}$,
the \emph{curvature}, satisfying
\[
 d \omega + \frac{1}{2} \left[\omega,\omega\right]
=\frac{1}{2} K \sigma \wedge \sigma.
\]

Let $V_0=\mathfrak{g}/\mathfrak{h}$.
Map $E \to FTM$ by $e \in E \mapsto \omega_e + \mathfrak{h} \in FT_{\pi(m)} M$.
This map is an $H$-structure. The 1-form
$\sigma$ is called the \emph{soldering form}.

Take the the obvious projection
$\mathfrak{g} \to \mathfrak{g}/\mathfrak{h}$.
Take any linear splitting of vector spaces
$s : \mathfrak{g}/\mathfrak{h} \to \mathfrak{g}$.
The function $t(a,b)=[s(a),s(b)]+K(a,b)+\mathfrak{h}$ (for
$a, b \in \mathfrak{g}/\mathfrak{h}=V_0$) is a torsion function
just when $s$ satisf\/ies
\[
\Ad(h) \left[s(A),s(B)\right]+\mathfrak{h}= \left[s\left(\Ad(h) A\right),s\left(\Ad(h) B\right)\right]+\mathfrak{h}
\]
for all $A,B \in \mathfrak{g}/\mathfrak{h}$.
(Warning: $t$ is \emph{not} necessarily the same as
the object which is usually called the
\emph{torsion} of the $G/H$-geometry. Moreover
it is not clear that there is always such a splitting
$s$ satisfying this complicated condition.)
The intrinsic torsion
$[t]$ is of course independent of the choice of splitting $s$.
The prolongation of the underlying f\/irst
order structure is $E \times_H \prol{\mathfrak{h}}$.
\end{Example}

\begin{Definition}
Suppose that $G$ is a complex Lie group with Lie algebra $\mathfrak{g}$
and $V_0$ is a complex $G$-module.
Recall the exact sequence
\[
\xymatrix{ 0 \ar[r] & \prol{\mathfrak{g}} \ar[r] & \mathfrak{g} \otimes V^*
  \ar[r]^{\delta} & V_0 \otimes \Lm{2}{V_0}^* \ar[r]^{[]}
& \Cohom{0,2}{\mathfrak{g}} \ar[r] &
  0.  }
\]
Map $\prol{\mathfrak{g}} \to \GL{V_1}$ by
\[
Q \mapsto Q', Q'(v,A) = (v,A+Qv).
\]
Suppose that $E \to FTM$ is a $G$-structure with torsion function $t$
and prolongation $\prol{E}$. Let $\prol{\pi} : \prol{E} \to E$ be the bundle map.
Let $FE$ be the $V_1$-valued frame bundle of $E$.
Map $\prol{E} \to FE$ by taking
$\gamma \mapsto \pi(\gamma)^* \sigma \oplus \gamma$.
Then $\prol{E} \to FE$ is a $\prol{\mathfrak{g}}$-structure on $E$,
called the \emph{prolongation} of $E$.

There is a right action of $Q \in \prol{\mathfrak{g}}$ given
by $r_Q U = U-Q \sigma$ (where $U \in \prol{E}$).
This action satisf\/ies
\begin{gather*}
  r_Q^* \sigma  = \sigma,\qquad
  r_Q^* \gamma  = \gamma - Q \sigma.
\end{gather*}
There is a natural action of $G$ on $\prol{E}$, commuting with the
bundle map $\prol{E} \to E$, which is (for $g \in G$):
\[
r_g U = \Ad_g^{-1} \left ( U \left ( r_g^{-1} \right )' \right ).
\]
Form the semidirect product
$G \rtimes \prol{\mathfrak{g}}$ with multiplication
\[
\left ( g_1, Q_1 \right ) \left ( g_2, Q_2 \right ) = \left ( g_1g_2,
  Q_1 + g_1 Q_2 \right )
\]
where $gQ$ means the element of $\Sym{2}{V_0}^* \otimes V_0$ def\/ined by
\[
gQ (u,v) = g \left ( Q \left( g^{-1} u, g^{-1} v \right ) \right ).
\]
We can write a point of $\prol{E}$ as $(u,U)$ where $u \in E$ and $U$ is a
pseudoconnection at $u$. The two group actions fuse together to an action of the
semidirect product:
\[
r_{(g,Q)} \left(u,U\right) = \left( g^{-1} u, \Ad_{g}^{-1} \left( U -
    Q u \pi' \right ) \left ( r_g^{-1} \right )' \right).
\]
This action makes $\prol{E} \to M$ into a principal right $G \rtimes
\prol{\mathfrak{g}}$-bundle. (It is \emph{not} a
$G \rtimes \prol{\mathfrak{g}}$-structure.) Under
the $G$-action
\begin{gather*}
  r_g^* \sigma   = g^{-1} \sigma, \qquad
  r_g^* \gamma   = \Ad_g^{-1} \gamma.
\end{gather*} Therefore under the $G \rtimes \prol{\mathfrak{g}}$-action
\[
r_{(g,Q)}^*
\begin{pmatrix}
  \sigma \\
  \gamma
\end{pmatrix}
=
\begin{pmatrix}
  g^{-1} \sigma \\
  \Ad_{g}^{-1} \left ( \gamma - Q \sigma \right )
\end{pmatrix}.
\]
\end{Definition}

\begin{Definition}\label{definition:SecondOrderStructure}
A \emph{second order structure} on a manifold $M$ with
f\/irst order structure $E \to FTM$ and torsion function $t$
and $G$-submodule $G_1 \subset \prol{\mathfrak{g}}$
is a principal $G_1$-bundle $E_1 \to E$
and $G \rtimes G_1$-equivariant
bundle map $E_1 \to \prol{E}$ over $M$.
(In this def\/inition, $\prol{E}$ is the prolongation of the
f\/irst order structure using $t$ as torsion function.)
\end{Definition}
\begin{Definition}
Suppose that $\phi_1 : E_1 \to E$ is a second order
structure over a $G$-structure $\phi : E \to M$.
We def\/ine $\Gromov{\phi}_1 : \prol{E} \to \prol{\mathfrak{g}}/G_1$
by $\Gromov{\phi}_1(\gamma) = Q+G_1$ if
$Q \in \prol{\mathfrak{g}}$ and $r_Q \gamma \in E_1$.
Much as before, $\Gromov{\phi}_1^{-1} G_1=E_1$.
\end{Definition}

We leave the reader to def\/ine higher order structures of all
orders by induction.

\begin{Proposition}\label{proposition:HigherOrder}
A holomorphic higher order structure
with discrete kernel extends holomorphically across a
subset of \cct
just when its underlying first order structure extends holomorphically
across that subset.
\end{Proposition}
\begin{proof}
It is enough (by induction) to prove the result for
second order structures.
Since the result is local, replace the complex manifold $M$ by a ball $B$. As usual,
because $E \to B \backslash S$ extends holomorphically,
say to a principal $G$-bundle $E' \to B$
we can replace $B$ by a smaller ball
to arrange that $E' \to B$ is holomorphically trivial, $E' = B \times G$,
and so $E = \left(B \backslash S\right) \times G$.
Therefore $\prol{E}=\left(B \backslash S\right) \times G \times \prol{\mathfrak{g}}$.
Over $E$ we have a $G_1$-structure $E_1 \to E$,
i.e. $E_1 \to \left(B \backslash S\right) \times G$.
So we have a map $\Gromov{\phi}_1 : \prol{E} \to \prol{\mathfrak{g}}/G_1$,
i.e. $\Gromov{\phi}_1 :  \left(B \backslash S\right) \times G \times
\prol{\mathfrak{g}} \to \prol{\mathfrak{g}}/G_1$.
The quotient $\prol{\mathfrak{g}}/G_1$ is a quotient
of vector spaces, so a vector space, and therefore a
\Thullen extension
target. Therefore for each choice of $g \in G$
and $Q \in \prol{\mathfrak{g}}$, this map
$\Gromov{\phi}_1$ extends to a holomorphic
map on $B$, and therefore
extends to a holomorphic map on
$B \times G \times \prol{\mathfrak{g}}=\prol{E}$.
We extend $E_1$ holomorphically to be $E_1'=\Gromov{\phi}_1^{-1}(G_1)$.
\end{proof}

\begin{Theorem}\label{theorem:HigherHartogs}
Let $M$ be a domain in a Stein manifold
and $\hat{M}$ the envelope of holomorphy of~$M$.
A holomorphic higher order structure
with trivial kernel extends holomorphically from~$M$
to~$\hat{M}$ just when its underlying first order
structure extends holomorphically
from~$M$ to~$\hat{M}$.
\end{Theorem}
\begin{proof}
It is enough (by induction) to prove the result for
second order structures.
Suppose that our higher order structure is $E_1 \to E \to M$,
in the notation of Def\/inition~\vref{definition:SecondOrderStructure}.
Suppose that $E$ extends to a bundle over $\hat{M}$,
which we denote by $E'$.
Consider the torsion function $t : E \to \Lm{2}{V_0}^* \otimes V_0$.
The torsion function is a section of the holomorphic
vector bundle $E \times_G \left( \Lm{2}{V_0}^* \otimes V_0 \right )$.
By Proposition~\vref{proposition:VectorBundleSections},
this section extends to a section of
$E' \times_G \left( \Lm{2}{V_0}^* \otimes V_0 \right )$,
and therefore the torsion extends to a holomorphic function
on $E'$.

By Lemma~\vref{lemma:SteinHasHolomorphicConnections},
there is a connection for $E' \to \hat{M}$, say~$\alpha$.
If $\alpha$ has torsion, say~$t_{\alpha}$,
and $E$ has torsion function $t$, then
replace $\alpha$ by $\alpha - \frac{1}{2} t_{\alpha} \sigma + \frac{1}{2} t \, \sigma$.
Every pseudoconnection 1-form~$U$ at any
point $u \in E$ has the form $U=\alpha + p \, \sigma$,
where $\sigma$ is the soldering 1-form and
$p \in \C{n*} \otimes \mathfrak{g}$. The torsion of $U$
is then $t+\delta p$, so we will need $p \in \prol{\mathfrak{g}}$.
In order that $U \in E_1$,
we need to pick $p$ to lie in some $G_1$-orbit
in $\prol{\mathfrak{g}}$. This orbit may vary
at dif\/ferent points of $E$. We can keep track
of this orbit as a function $f : E \to \prol{\mathfrak{g}}/G_1$.
This function is a section of
$E \times_G \left( \prol{\mathfrak{g}}/G_1 \right)$
so once again extends uniquely holomorphically
to a function $f : E' \to \prol{\mathfrak{g}}/G_1$.
We then produce a bundle $E_1' \to E'$
consisting precisely of the pseudoconnections
$U=\alpha+p \, \sigma$ for which $\delta p=0$
and $p + G_1 = f$.
\end{proof}

\subsection{Parabolic geometries and higher order structures}

\begin{Example}\label{example:ProjectiveConnectionHartogsAndThullen}
Consider $G=\GL{n,\C{}}$. Then $\prol{\mathfrak{g}}=V_0 \otimes \Sym{2}{V_0}^*$.
There is an obvious $G$-submodule $G_1=V_0^*$ included by taking $a \in V_0^*$
to the operator
$\left(v_1,v_2\right) \mapsto a\left(v_1\right)v_2+a\left(v_2\right)v_1$.
A \emph{normal projective connection} on a manifold $M$
is a second order
structure with f\/irst order
structure $E=FTM$, torsion function $t=0$,
and $G$-submodule $G_1 \subset \prol{\mathfrak{g}}$.
See Borel \cite{Borel:2001}
or Molzon and Mortensen \cite{Molzon/Mortensen:1996}
for more information.
Normal projective connections extend across
subsets of \cct, and extend from
domains in Stein manifolds to envelopes
of holomorphy, because the
underlying f\/irst order structure is just the frame bundle.
In fact, normal projective
connections turn out to be equivalent to Cartan connections
modelled on projective space (again see Borel \cite{Borel:2001}), so we
will also solve the \Thullen extension problem for
normal projective connections in
Theorem~\vref{theorem:ParabolicGeometriesExtend}.
\end{Example}

\begin{Example}\label{example:ContactProjectiveConnectionHartogsAndThullen}
A \emph{normal contact projective connection} is a similar
sort of higher order structure, to complicated
to explain in detail; see Fox~\cite{Fox:2005}.
Its underlying f\/irst order structure is a contact
structure, so that we have solved the \Thullen
and Hartogs extension problems for normal contact projective
connections by combining
Theorem~\vref{theorem:ContactStructuresExtend} (to
extend the contact structure) with
Proposition~\vref{proposition:HigherOrder} (to
extend the second order structure across
subsets of \cct) or Theorem~\vref{theorem:HigherHartogs}
(to extend from a domain in Stein manifold
to the envelope of holomorphy).
Normal contact projective
connections turn out to be equivalent to Cartan connections
modelled on the projectivized cotangent bundle of
projective space (again see Fox~\cite{Fox:2005}), so we
will also solve the \Thullen extension problem for
normal contact projective connections in
Theorem~\vref{theorem:ParabolicGeometriesExtend}.
\end{Example}

\section{Extending f\/lat Cartan geometries}

\subsection{Extending by development}

\begin{Definition}
Consider a $G/H$-geometry $\pi : E \to M$,
and let $\mathfrak{h} \subset \mathfrak{g}$
be the Lie algebras of $H \subset G$.
The 1-form $\sigma = \omega + \mathfrak{h} \in \nForms{1}{E} \otimes
\left(\mathfrak{g}/\mathfrak{h}\right)$ is called
the \emph{soldering form} of the Cartan geometry.
The curvature can be written as
$d \omega + \frac{1}{2}\left[\omega,\omega\right]
= \frac{1}{2} K \sigma \wedge \sigma$
for a unique function
$K : E \to \mathfrak{g} \otimes \Lm{2}{\mathfrak{g}/\mathfrak{h}}^*$,
called the \emph{curvature function}. The curvature function
transforms under $H$-action according to the obvious
equivariance, so is a section of
\[
 E \times_H \left(
\mathfrak{g} \otimes \Lm{2}{\mathfrak{g}/\mathfrak{h}}^*
\right),
\]
which we will also call the \emph{curvature}. We will
call a Cartan geometry \emph{flat} if its curvature
vanishes.
\end{Definition}

\begin{Lemma}
A Cartan geometry is locally isomorphic to its model just when
it is flat.
\end{Lemma}
This lemma is well known
(see Sharpe \cite[Theorem~5.1, p.~212]{Sharpe:1997})
and doesn't require complex analyticity.
\begin{proof}
Take $E \to M$ the Cartan geometry, with Cartan
connection $\omega$. Apply the Frobenius theorem
to the Pfaf\/f\/ian system $\omega-g^{-1} \, dg=0$ on $E \times G$.
The integral manifolds of this Pfaf\/f\/ian system are locally graphs of local isomorphisms.
\end{proof}

\begin{Lemma}
Suppose that $E \to M$ is a flat $G/H$-geometry.
Let $H$ act on $E \times G$ by the right
action $(e,g)h=\left(eh,gh\right)$.
Take any connected integral manifold $Z$ of the Pfaffian system
$\omega-g^{-1} \, dg=0$ on $E \times G$.
Let $ZH$ be the union of all $H$-orbits through points of $Z$.
Then $ZH$ is an $H$-equivariant covering space of $E$,
and the total space of a Cartan geometry on a covering space of $M$.
\end{Lemma}
\begin{proof}
We provide a sketch; see McKay \cite{McKay:2007}
for details. Above each curve in $E$,
the Pfaf\/f\/ian system is an ordinary dif\/ferential equation
of Lie type, so has global solution. Therefore each integral
manifold $Z$ is a covering space of a path component of $E$.
The $H$-orbits are Cauchy characteristics of the Pfaf\/f\/ian system.
Therefore the group of path components of $H$ acts permuting
integral manifolds. The union $ZH$ of all of these $H$-orbits
is therefore a union over path components of $H$, say
$ZH=Z\pi_0(H)$, a discrete union of distinct connected integral
manifolds, so an integral manifold. Moreover, $ZH$ is acted
on by $H$ freely and properly, because~$H$ acts freely
and properly on $E \times G$. Therefore $E' = ZH \to M' = ZH/H$
is a Cartan geometry with Cartan connection $\omega$.
\end{proof}

\begin{Proposition}\label{prop:WeakFlatCartan}
Suppose that $G/H$ is a complex homogeneous
space. Then local biholomorphisms to $G/H$ extend
across subsets of \cct if and only if
flat holomorphic $G/H$-geometries extend across
such subsets.
\end{Proposition}
\begin{proof}
Suppose that $G/H$ is a \Thullen extension target for local biholomorphisms.
We can replace $M$ by any open neighborhood of a point
$s \in S$, so we can assume that $M$ is a ball $B$.
Then $B \backslash S$ is simply connected.
If $Z$ is any integral manifold of the Pfaf\/f\/ian
system $\omega-g^{-1} \, dg=0$ on $E \times G$, then $ZH$ is
the graph of a local isomorphism $f: B  \backslash S \to G/H$.
Because $G/H$ is a \Thullen extension target
for local biholomorphisms, we can
extend the map $f$. The pullback
of the bundle $G \to G/H$ holomorphically extends
the bundle $E$. By Theorem~\vref{theorem:CartanGeometry}, the Cartan
geometry extends holomorphically.

Suppose that some local biholomorphism $f : M \backslash S \to G/H$
does not extend to $M$, with $S \subset M$ a subset of \cct.
Then the pullback $E = f^* G$ of the standard
f\/lat $G/H$-geometry is a $G/H$-geometry on $M \backslash S$.
To be precise, $E$ is the set of pairs
$\left(m_0,g_0\right) \in \left(M \backslash S\right) \times G$
so that $f\left(m_0\right)=g_0H \in G/H$.
Suppose that this $G/H$-geometry extends to a $G/H$-geometry $E' \to M$.
The curvature vanishes on $E$, a dense open set in $E'$,
so the $G/H$-geometry on $E'$
is f\/lat. Pick a maximal integral manifold $Z$ of the
Pfaf\/f\/ian system $\omega-g^{-1} \, dg=0$ on $E' \times G$
passing through a point of the f\/iber $E'_s$.
This integral manifold $Z$ is uniquely determined
up to left $G$-action on $G$ and right $H$-action on $E'$.
After perhaps translating
$Z$ by right $H$-action, we can arrange that $Z$ passes
through a point of the form $\left(m_0,g_0,g_0\right) \in E \times G$.
This leaf~$Z$ is then locally the graph of the
local isomorphism $\left(m_0,g_0\right) \in E \mapsto g_0 \in G$
of $G/H$-geometries, and extends this isomorphism to
a neighborhood of $E'_s$ in $E'$. Therefore $Z$ is the
graph of a local isomorphism of $G/H$-geometries
taking an open neighborhood of $s \in M$ to an open
set in~$G/H$, extending $f$.
\end{proof}

\subsection{Examples of inextensible f\/lat Cartan geometries}

\begin{Lemma}[Sharpe \protect{\cite[Theorem 3.15, p.~188]{Sharpe:1997}}]\label{lemma:TgtBundle}
If $\pi : E \to M$ is any $G/H$-geometry
then the Cartan connection of $E$ maps
\[
\xymatrix{%
0 \ar[r] & \ker \pi'(e) \ar[r] \ar[d] & T_e E \ar[r] \ar[d] & T_m M \ar[r] \ar[d] & 0 \\
0 \ar[r] & \mathfrak{h} \ar[r] & \mathfrak{g} \ar[r] &
\mathfrak{g}/\mathfrak{h} \ar[r] & 0 }
\]
for any points $m \in M$ and $e \in E_m$; thus
\[
TM=E \times_H \left(\mathfrak{g}/\mathfrak{h}\right).
\]
Under this identification, holomorphic vector fields are identified
with $H$-equivariant holomorphic functions $E \to \mathfrak{g}/\mathfrak{h}$.
\end{Lemma}

\begin{Example}
Reconsider Example~\vref{example:proj}.
Take $G \subset \GL{n,\C{}}$ any Lie group
acting transitively on $\C{n} \backslash 0$,
and let $H$ be the stabilizer
of some vector $v_0 \in \C{n} \backslash 0$.
Obviously $G/H=\C{n} \backslash 0$.
For example, $G=\Symp{2n,\C{}}$ on $\C{2n}$,
or $G=G_2$ on $\C{7}$.
Then the standard f\/lat $G/H$-geometry,
on $G/H=\C{n} \backslash 0$ does not extend across
the puncture at $0$.
For any $G/H$-geometry $E \to M$,
the constant function $v_0$ on $E$
is identif\/ied with a nowhere vanishing vector f\/ield
by Lemma~\ref{lemma:TgtBundle}.
On the model $G/H=\C{n} \backslash 0$,
this vector f\/ield is the Euler vector
f\/ield $Z(z)=z$. Clearly the Euler
vector f\/ield cannot extend to $\C{n}$
holomorphically without vanishing at the
origin. If the $G/H$-structure were
to extend holomorphically to the origin,
then the associated nowhere vanishing
vector f\/ield would also extend holomorphically,
and remain nowhere vanishing.
\end{Example}

\section{Extending Cartan geometries}

\begin{Theorem}
 Suppose that $H \subset G$ is closed complex
Lie subgroup of a complex Lie group. If
holomorphic maps to $G/H$ extend across subsets of \cct,
then $G/H$-geometries extend
holomorphically across such subsets.
\end{Theorem}
This solves the \Thullen extension problem for various
Cartan geometries; the
analogous Hartogs extension problem is unsolved.
\begin{proof}
As usual we can assume that our $G/H$-geometry
is on a ball minus a subset of \cct, say
$E \to B \backslash S$.
Fix a point $z_0 \in B \backslash S$.
Consider the complex af\/f\/ine lines through $z_0$.
Away from $z_0$, each point of $B$ lies on a
unique disk from this family.

Pick a point $e_0 \in E_{z_0}$ in the f\/iber above $z_0$,
and develop each disk so that
the frame $e_0$ is carried to the frame $1 \in G$.
The development of each disk is well def\/ined,
by Lemma~\vref{lemma:developCurves},
except possibly for development along
the disks through $z_0$ which hit points of $S$.

The disk $D$ through $z_0$ and some point of $S$
has to be punctured at all points of $D \cap S$
before we develop,
because the Cartan geometry is not def\/ined
on $S$. The development
is def\/ined on the universal covering space
of the punctured disk. However, by continuity
of solutions of ordinary dif\/ferential equations,
the monodromy must be the limit of the monodromies
of nearby punctured disks. The nearby disks have
no punctures, so the monodromy is trivial:
the development is a~holomorphic map on each disk.

These developments f\/it
together into a single holomorphic map $\phi_1 : B \backslash S \to G/H$.
The development of each disk, say $D$, yields an isomorphism
$\Phi_D : \left.E\right|_D \to \left. \phi_1^*G \right|_D$.
This map clearly extends to a bundle isomorphism
$\Phi : E \to \phi_1^*G$ above $B \backslash S$. Since
the map $\phi_1$ extends across $S$, so
does the bundle $\phi_1^*G$, and therefore so does $E$.
\end{proof}

\begin{Example}\label{example:reductiveAlgebraicCartan}
If $G/H$ is a reductive homogeneous space, then $G/H$-geometries
extend across subsets of \cct.
\end{Example}

\begin{Example}\label{example:RiemannCartanSurfaces}
Among Cartan geometries on surfaces, this theorem applies to
\begin{itemize}\itemsep=0pt
\item $\C{2}$-geometries, for any of
the various complex Lie groups acting transitively
on $\C{2}$ and
\item $\Proj{1} \times \Proj{1} \backslash \text{diagonal}$-geometries;
see Example~\vref{example:web}.
\end{itemize}
\end{Example}

\begin{Example}
Beloshapka, Ezhov and Schmalz \cite{Beloshapka/Ezhov/Schmalz:2005}
study a particular class of {CR}-manifolds, those of dimension~4, CR-dimension~1, CR-codimension~2, with Engel CR-distribution.
They associate to each such CR-structure a canonical Cartan
geometry modelled on a homogeneous space $G/H$. If we
complexify this type of geometry, the resulting
complex analytic geometric structures are Cartan
geometries modelled on the complexif\/ication of $G/H$.
This complexif\/ication turns out to be an af\/f\/ine
space, and therefore these complex analytic
geometries extend across subsets
of \cct.
\end{Example}

\begin{Theorem}\label{theorem:ParabolicGeometriesExtend}
Parabolic geometries extend across subsets of \cct.
\end{Theorem}
\begin{proof}
 Pick a complex manifold $M$, a subset $S \subset M$ of \cct,
and a~parabolic geometry $E \to M \backslash S$
modelled on some rational homogeneous variety $G/P$.
We only need to extend $E$ locally, so we can assume that $M$
admits a holomorphic volume form. By McKay \cite[Theorem~2,  p.~2]{McKay:2008},
because $M \backslash S$ has a holomorphic volume form,
$E \to M \backslash S$ admits a
holomorphic reduction $E_0 \subset E$ of structure
group to a reductive algebraic group $P_0 \subset P$,
and~$E_0$ has a holomorphic connection.
By Proposition~\vref{proposition:connection},
$E_0$ extends to a holomorphic bundle $E_0' \to M$ with holomorphic
connection. Then $E$ extends to $E'=E_0' \times_{P_0} P$.
\end{proof}

\begin{Remark}
Hong \cite{Hong:2000} proved a related result.
She assumes
the existence of a family
of minimal rational curves on a
Fano manifold, satisfying
some complicated conditions on their
tangent lines at a generic point.
These conditions turn out to encode
a f\/lat parabolic geometry outside a
subset of the Fano manifold of \cct. She demonstrates
that the parabolic geometry extends holomorphically along
the curves, to become holomorphic at
every point of the manifold.
It is possible to prove Hong's theorem
by using the theorem above,
but one would f\/irst need
results of Ochiai~\cite{Ochiai:1970}
to demonstrate that the geometry
encoded by the tangents of the
minimal rational curves is actually a
parabolic geometry. We will leave this
to the reader. Mok~\cite{Mok:2002} made use of this
result, and of the result of Hwang and Mok
on codimension 2 extension of reductive $G$-structures,
in his study of the contact parabolic
geometry of $G_2$.
\end{Remark}

\begin{Example}
This theorem provides another proof that local biholomorphisms
to rational homogeneous varieties extend across subsets of \cct.
\end{Example}

\begin{Example}
 This theorem is an application of the method of harmless reduction
(see Section~\vref{section:HarmlessReductions}).
\end{Example}

\begin{Example}\label{example:ODEs}
There are (at least) two natural geometric structures associated to
a 2nd order scalar ordinary dif\/ferential equation. We
contrast their extension theory.

Recall the complex homogeneous surface $\OO{n}$
(see Example~\vref{example:homogeneousSurfaces}).
It turns out (see Dunajski and Tod \cite{Dunajski/Tod:2005}) that
every complex analytic scalar ordinary
dif\/ferential equation of order $n+1$
\[
 \frac{d^{n+1}w}{dz^{n+1}} = f\left(z,w,\frac{dw}{dz}, \dots, \frac{d^n w}{dz^n}\right)
\]
imposes an $\OO{n}$-geometry on the open subset
of the $(z,w)$-plane where the equation is def\/ined,
invariant under f\/iber-preserving transformations,
i.e.\ transformations preserving the solutions
of the dif\/ferential equation and the vertical
lines $z=\text{const}$. Every $\OO{n}$-geometry
with vanishing torsion arises in this fashion;
again see Dunajski and Tod \cite{Dunajski/Tod:2005} for a def\/inition of
torsion and a~proof.

An $\OO{n}$-geometry on a surface has as its underlying
f\/irst order structure a foliation with af\/f\/ine structure
on its leaves. This foliation comes from the
invariant foliation (indeed line bundle mapping)
$\OO{n} \to \Proj{1}$. The af\/f\/ine structure of the f\/ibers
endows the leaves of any $\OO{n}$-geometry with af\/f\/ine structures.
We have already seen an example
of an $\OO{n}$-geometry
in Example~\vref{example:homogeneousSurfaces}:
map $\C{2} \backslash 0 \to \OO{n}$ by taking
a point $z \ne 0$ to the homogeneous polynomial $p_z$
of degree $n$ on the
span of $z$ which takes the value $1$ on $z$.
The foliation is the foliation of $\C{2} \backslash 0$
by lines through $0$, so doesn't extend across
the puncture at $0$.

Let $A$ be any $2 \times 2$ invertible complex
matrix which has all
of its eigenvalues inside the unit disk.
Consider the \emph{Hopf surface}
$S=\left(\C{2} \backslash 0\right)/(z \sim Az)$,
which is a smooth compact complex surface.
The $\OO{n}$-geometry on $\C{2} \backslash 0$
descends to the Hopf surface $S$.
There are inf\/initely
many holomorphic $\OO{n}$-geometries on certain Hopf
surfaces; the author hopes to classify
these in a~sequel to this paper.

A dif\/ferent approach: every 2nd order scalar ordinary dif\/ferential
equation
\[
 \frac{d^{2}w}{dz^{2}} = f\left(z,w,\frac{dw}{dz}\right)
\]
gives rise to a contact 3-fold, with coordinates say $z,w,p$,
and contact structure $dw-p \, dz=0$. This 3-fold has two
nowhere tangent Legendre f\/ibrations: (1) $dw-p \, dz=dp-f \, dz=0$ and
(2)~$dw=dz=0$. This contact structure and
pair of Legendre f\/ibrations is invariant under point
transformations of the dif\/ferential equation.
Conversely, we def\/ine a \emph{path geometry}
to be any pair of nowhere tangent holomorphic Legendre f\/ibrations
in a holomorphic contact 3-fold. Path geometries, being parabolic
geometries modelled on $\Proj{}T\Proj{2}$ (see McKay~\cite{McKay:2006}
for proof), extend across subsets of \cct.

A similar story occurs in 3rd order: a 3rd order scalar ordinary
dif\/ferential equation is equivalent to a parabolic geometry
on an appropriate manifold, and this equivalence is
invariant under contact transformations; see
Sato and Yoshikawa \cite{Sato/Yoshikawa:1998} for proof.
The induced $\OO{2}$-geometry of a 3rd order scalar ODE
might not extend across a subset of \cct, but the
contact invariant parabolic geometry does extend
across such subsets.
\end{Example}

\begin{Remark}
This striking dif\/ference between f\/iber-preserving
and point transformations is probably quite important
in understanding 2nd order ODEs, and probably the most
important observation in this paper.
\end{Remark}

\begin{Remark}
\v{C}ap and Schichl \cite{Cap/Schichl:2000}
proved that any parabolic geometry
satisfying a
a mild (but complicated) hypothesis is
a local product
structure, with each local factor
bearing a~normal projective connection,
normal contact projective connection,
or a f\/irst order structure consisting
of various subbundles of the tangent
bundle, equipped with a reduction
of structure group to a~reductive
algebraic group. Suppose that we
have a holomorphic parabolic geometry
on a~domain in a Stein manifold,
which satisf\/ies the mild hypothesis
of \v{C}ap and Schichl.
The local product
structure is a f\/irst order structure
with reductive structure group, so
extends to the envelope of holomorphy
(see Example \vref{example:AlmostProduct}).
The various subbundles of the tangent
bundle extend, except perhaps on
a subset of \cct
(see Example~\vref{example:PlaneFields}).
The normal contact projective connections
and normal projective connections
will extend by arguments which
only slightly generalize those
of Examples~\ref{example:ProjectiveConnectionHartogsAndThullen}
and~\vref{example:ContactProjectiveConnectionHartogsAndThullen}.
Parabolic geometries extend
across subsets of \cct,
so we can work modulo such subsets.

In order that extension of the parabolic geometry
be parabolic,
the various plane f\/ields need to have
constant symbol algebra, up
to complex codimension 2;
see \v{C}ap and Schichl \cite{Cap/Schichl:2000} for
more information. This is
not dif\/f\/icult to prove, since
the symbol algebra changes type
only on an analytic set.
The mild hypothesis
of \v{C}ap and Schichl will therefore
be  suf\/f\/icient
to solve the Thullen and Hartogs
extension problems. We hope to solve
both extension problems
for all holomorphic parabolic geometries
(without even a mild hypothesis)
in a sequel to this paper.
\end{Remark}

\section{Extension of local isomorphisms}

A related question: for which holomorphic geometric structures
do local isomorphisms
extend holomorphically across subsets of \cct?
\begin{Example}
We produce an example of a holomorphic local isomorphism
of $G$-structures that fails to extend across a point,
even though the $G$-structures are f\/lat,
the isomorphism maps to a compact complex homogeneous
space, and the homogeneous space $\GL{n,\C{}}/G$ is a Hartogs
extension target.

Let $G$ be the set of $n \times n$ matrices of the form
$2^k I$ for $k \in \mathbb{Z}$. A $G$-structure
is a choice of holomorphic framing of a complex manifold,
up to rescaling by factors of $2$.
The quotient $\GL{2,\C{}}/G$ is covered by $\GL{2,\C{}}$, which
is a \Thullen and Hartogs extension
target, and therefore $\GL{2,\C{}}/G$ is also
a \Thullen and Hartogs extension target. Therefore $G$-structures extend
holomorphically across subsets of \cct.

Let $S$ be the Hopf manifold $\left(\C{2} \backslash 0\right)/(z \sim 2z)$,
and take the map $f : \C{2} \backslash 0 \to S$ taking
each point $z$ to its equivalence class $[z] \in S$.
All linear automorphisms of $\C{2}$ commute with
$z \mapsto 2z$, so they descend to biholomorphisms
of the Hopf surface. These linear automorphisms act transitively
on $\C{2} \backslash 0$, so act transitively on the
Hopf surface.

Any local biholomorphism $f : M \to N$ between complex manifolds
determines a map $f_1 : FTM
\to FTN$ by $Ff(m,u)=u \circ f'(m)^{-1}$ for $m \in M$
and $u \in FT_m M.$ Since
$FT\left(\C{2} \backslash 0\right)
=\left(\C{2} \backslash 0\right) \times \GL{2,\C{}}$
is a trivial bundle, we can compose with the obvious global
section of that bundle, to obtain a map $f : \C{2} \backslash 0 \to FTS$,
a $G$-structure on $S$. We can map the standard f\/lat
$G$-structure on $\C{2} \backslash 0$ to the given
$G$-structure on $S$, also called $f$,
by
\[
 f\left(z,2^k I \right)= f(z)2^kI.
\]
This map is a local isomorphism of $G$-structures, taking
the standard f\/lat $G$-structure to the $G$-structure on the
Hopf surface, $\C{2} \backslash 0 \to S$. The map doesn't extend
holomorphically to $\C{2}$. Clearly we can generalize to
Hopf manifolds of all dimensions.
\end{Example}

\section{Rigidity}

Many of the theorems in this paper have
obvious generalizations to families
of f\/irst order structures, higher
order structures and Cartan geometries.
We will provide some examples to
expose the general pattern.

We will think of a holomorphic submersion
$\pi : M^{n+p} \to Z^p$ as a $p$-parameter family
of complex manifolds.

\begin{Definition}
The \emph{vertical bundle} $V=V_{\pi}$ of such
a holomorphic submersion $\pi : M \to Z$
is $V=\ker \pi' \subset TM$.
The \emph{vertical frame bundle} is the
$\GL{n,\C{}}$-bundle $FV \to M$.
Suppose that $\rho : G \to \GL{n,\C{}}$
is a complex representation of a complex
Lie group $G$.
A \emph{holomorphic family of $G$-structures}
on a holomorphic submersion $\pi : M^{n+p} \to Z^p$
is a holomorphic principal right $G$-bundle $E \to M$
and a $G$-equivariant holomorphic bundle map $E \to FV$.

If $p : E \to M$ is a holomorphic principal bundle
and $\pi : M \to Z$ is a holomorphic submersion,
let $VE \to E$ be the bundle
$VE = \ker \pi p \subset TE$.
Suppose that $G/H$ is a complex homogeneous space.
A \emph{family of Cartan geometries} modelled
on $G/H$ on a holomorphic submersion $\pi : M \to Z$
is a holomorphic principal right
$H$-bundle $E \to M$ and a section
$\omega$ of $VE \otimes \mathfrak{g}$,
called the \emph{Cartan connection},
satisfying the following conditions:
\begin{enumerate}\itemsep=0pt
\item
Denote the right action of $g \in G$ on $e \in E$ by $r_g e=eg$.
The Cartan connection transforms in the adjoint representation:
\[
r_g^* \omega = \Ad_g^{-1} \omega.
\]
\item
$\omega_e : V_e E \to \mathfrak{g}$ is a linear isomorphism at each point
$e \in E$.
\item
For each $A \in \mathfrak{g}$, def\/ine a vector f\/ield $\vec{A}$ on $E$,
tangent to the f\/ibers of $\pi p : E \to Z$, by
the equation $\vec{A} \hook \omega = A$. Then the vector f\/ields
$\vec{A}$ for $A \in \mathfrak{h}$ must generate the right $H$-action:
\[
\vec{A} = \left. \frac{d}{dt} r_{e^{tA}} \right|_{t=0}.
\]
\end{enumerate}
If $G/P$ is a rational homogeneous variety,
then a family of Cartan geometries modelled on
$G/P$ is called a family of \emph{parabolic geometries}.
\end{Definition}

Recall the foliation of $\C{n} \setminus 0$
by radial lines.
Clearly foliations, and even holomorphic
submersions, do not generally extend to
envelopes of holomorphy. Therefore
we will pose the Hartogs extension
problem only for holomorphic submersions
which are already assumed to extend.

\begin{Theorem}\label{theorem:ReductiveGstructuresHartogsExtensionFamilies}
Suppose that $G \subset \GL{n,\C{}}$ is a reductive
algebraic group.
Suppose that $M^{n+p}$ and $Z^p$ are both domains in Stein manifolds,
with envelopes of holomorphy $\hat{M}$ and $\hat{Z}$
respectively.
Suppose that $\pi : \hat{M} \to \hat{Z}$
is a holomorphic submersion.
Then every holomorphic family of $G$-structures on $\pi : M \to Z$ extends
uniquely to a holomorphic family of $G$-structures on $\pi : \hat{M} \to \hat{Z}$.
\end{Theorem}
\begin{proof}
Let $VM$ be the vertical bundle $\ker \pi'$ of $\pi : M \to Z$
and $V \hat{M}$ be vertical bundle $\ker \pi'$ of $\pi : \hat{M} \to \hat{Z}$.
A $G$-structure on $M$ is equivalent to a section of
$FVM/G \subset FV\hat{M}/G$.
The total space of $FV\hat{M}/G$ is
a Hartogs extension target by
Corollary~\vref{corollary:ExtendReductiveReduction}.
Therefore if $s : M \to FVM/G$ is a $G$-structure,
then $s$ extends uniquely to a holomorphic map
$s : \hat{M} \to FV\hat{M}/G$.
Let $p : FV\hat{M}/G \to \hat{M}$ denote
the bundle map. Then $p s$ is the identity
on $M$, and therefore by analytic continuation
is the identity on $\hat{M}$. So the
extension is also a holomorphic family of $G$-structures.
\end{proof}

The reader should compare the proof to that of
Theorem~\vref{theorem:ReductiveGstructuresHartogsExtension};
it is a small change of notation. We give one
more example:

\begin{Theorem}\label{theorem:ParabolicGeometriesExtendFamilies}
Families of parabolic geometries extend across subsets of \cct.
\end{Theorem}

First we will need a lemma from representation theory:

\begin{Definition}
Suppose that $P \subset G$ is a parabolic
subgroup of a complex semisimple Lie group.
There is a Cartan subgroup of $G$, say $H$,
which lies in $P$. Fixing a choice of
Cartan subgroup, we induce a choice of
root system for $G$. Divide up the roots
of $G$ into those whose root spaces lie in
$\mathfrak{p}$ and the rest. The former
we will call $P$-roots, and their root
spaces we will call $P$-root spaces. Among
the $P$-roots $\alpha$, there are
those for which $-\alpha$ is also
a $P$-root; call these $M$-roots.
The remaining $P$-roots we call
$N$-roots. The roots which are
not $P$-roots we call $N^{-}$-roots.
\end{Definition}

\begin{Lemma}[Langlands]\label{lemma:Langlands}
Suppose that $P \subset G$ is a parabolic
subgroup of a complex semisimple Lie group.
Then $P$ admits a decomposition
$P=MAN$ into connected complex
Lie subgroups, where $M$ is a complex
semisimple Lie group, $A$ is a
abelian linear algebraic Lie group, $MA$ is a maximal
reductive algebraic subgroup of $P$,
and $N$ is a unipotent subgroup.
At most a finite subgroup of $MA$ acts trivially
on $\mathfrak{g}/\mathfrak{p}$.
\end{Lemma}
For proof see Knapp \cite[p.~478]{Knapp:2002}.

Now we prove Theorem~\ref{theorem:ParabolicGeometriesExtendFamilies}.
\begin{proof}
It is clearly that solutions of Hartogs
or \Thullen extension problems are
unique when they exist, so we can
treat this problem as a local problem.
Pick a complex manifold $M$, a~subset $S \subset M$ of \cct,
a holomorphic submersion $\pi : M \to Z$, with vertical
bundle $VM \to M$,
and a family of parabolic geometries $E \to M \backslash S$
modelled on some rational homogeneous variety $G/P$.
In terms of the Langlands decomposition, $P=MAN$,
the maximal reductive subgroup of $P$ is $MA$.
We only need to extend $E$ locally, so we can assume that
$\det V^*M$
admits a nowhere vanishing holomorphic section. By McKay \cite[Theorem 2, p.~2]{McKay:2008},
$E \to M \backslash S$ admits a
holomorphic reduction $E_0 \subset E$ of structure
group to the maximal reductive algebraic group $MA \subset P$.
By Lemma~\vref{lemma:Langlands},
the group $MA$ acts on $\mathfrak{g}/\mathfrak{p}$
with only some f\/inite subgroup $K \subset MA$
acting trivially.
Therefore the underlying family of f\/irst order structures
is exactly an embedding $E_0/K \to FV\left(M \setminus S\right)$.

Since the problem is local, we can assume
that $M$ is a ball in $\C{n+p}$. But then
$M \setminus S$ has envelope of holomorphy
$M$.
By Theorem~\vref{theorem:ReductiveGstructuresHartogsExtensionFamilies},
the underlying f\/irst order structure
extends to a bundle mapping $E_0'' \to FVM$ for a
holomorphic principal right $MA/K$-bundle $E_0'' \to M$.
The covering bundle $E_0 \to E_0'$ extends to a covering
bundle $E_0' \to E_0''$, say, by Lemma~\vref{lemma:coveringSpacesExtendCCT}.
Let $E' = E_0' \times_{MA} P$.
By Proposition~\vref{proposition:VectorBundleSections}, we can extend
the holomorphic section $\omega$ of
the bundle $\left(VE \otimes \mathfrak{g}\right)^H \to M \setminus S$
to a holomorphic section of
$\left(VE' \otimes \mathfrak{g}\right)^H \to M$.

This section is thus an $H$-invariant
section of $VE \otimes \mathfrak{g}$, the f\/irst property
of a Cartan connection for a family of
Cartan geometries.

Let $H$ be the set of points $e \in E'$ at
which $\omega_e : V_e E' \to \mathfrak{g}$
is \emph{not} a linear isomorphism. Clearly
$H$ is a hypersurface, given by the one equation
$\det \omega_e=0$. Moreover, this hypersurface
is $P$-invariant, so projects to a hypersurface
in $\hat{M}$. This hypersurface doesn't intersect
$M$, so is empty by Lemma~\vref{lemma:hypersurfaceIntersections}.
Therefore $\omega$ satisf\/ies the second property
of a Cartan connection.

Over $M$, $\omega$ satisf\/ies $\vec{A} \hook \omega = A$,
for any $A \in \mathfrak{h}$. By analytic continuation,
this must also hold over $\hat{M}$, the third and f\/inal property
of a Cartan connection.
\end{proof}

\begin{Example}
\begin{Proposition}
Pick a rational homogeneous variety $G/P$.
Pick a holomorphic fibration $\pi : M^{n+p} \to Z^p$
whose fibers are biholomorphic to $G/P$ away from
a set of \cct in $Z$.
There is a holomorphic principal
$G$-bundle $B \to Z$ so that
this fibration is a~holomorphic fiber bundle $M = B/P \to Z$.
\end{Proposition}
\begin{proof}
Suppose that $\pi : M^{n+p} \to Z^p$ is a holomorphic
f\/ibration. Let $S_Z \subset Z$ be a subset
of \cct, and let $S_M=\pi^{-1} S_Z$.
Suppose that each f\/iber $M_z \subset M$
is biholomorphic to $G/P$, as long as $z$
is not a point of $S_Z$. A rational
homogeneous variety admits a unique
$G/P$-geometry: the standard f\/lat one;
see McKay \cite[Corollary~7, p.~16]{McKay:2006}.
Therefore $\pi : M \setminus S_M \to Z \setminus S_Z$
admits a unique holomorphic
family of $G/P$-geometries.
By Theorem~\vref{theorem:ParabolicGeometriesExtendFamilies},
this family of parabolic geometries
extends holomorphically
to a family of parabolic geometries on
$M$. Therefore each of the f\/ibers $M_s$ for $s$ in $S_Z$
also bears a parabolic geometry. By continuity
of the curvature, the parabolic geometry on each f\/iber
$M_s$ is f\/lat.
Every holomorphic f\/ibration is a
f\/iber bundle of smooth real manifolds.
Since $G/P$ is compact and simply connected,
each f\/iber $M_s$ is also compact and simply connected.
Therefore $M_s$ has a~holomorphic developing
map to $G/P$, a local biholomorphism. Because
the f\/iber $M_s$ is compact, the developing
map is a covering map. Because $G/P$ is simply
connected, the developing map is a biholomorphism.
Therefore all f\/ibers are biholomorphic
to $G/P$, and the biholomorphisms preserve
the parabolic geometry, so preserve the
$G$-action. Local triviality is
clear by using the developing map on a local
section.
\end{proof}
This is a mild form
of rigidity for $G/P$, much weaker than results of Hong~\cite{Hong:2000}
and of Hwang and Mok~\cite{Hwang/Mok:1998,Hwang/Mok:2005}.
\end{Example}

\section{Conclusion}

Holomorphic geometric structures
which occur naturally
apparently almost always extend across subsets
of \cct, and extend from a domain
in a Stein manifold to its envelope of holomorphy.
I found the idea behind the Hartog's lemma
for Cartan geometries in Gunning \cite[p.~126]{Gunning:1978};
the same idea appears
for $G$-structures in Hwang and Mok~\cite{Hwang/Mok:1998}.
Similarly, you might hope to extend
Cartan geometries or $G$-structures across real or complex analytic
submanifolds and analytic subsets using well known
extension techniques; see Siu~\cite{Siu:1974} for an introduction
to these techniques. For example,
restrict holomorphic
geometric structures to the boundary of a pseudoconvex
domain, to obtain geometric structures on the boundary.
It seems natural to ask which geometric structures on the boundary
arise in this way.

The \emph{generalized Cartan geometries}
of Alekseevsky and Michor~\cite{Alekseevsky/Michor:1995}
have the same def\/inition
as Cartan geometries (see Def\/inition~\vref{def:CartanConnection}),
except that they omit condition~2. Examples arise naturally
from maps to homogeneous spaces, and also from compactif\/ication
problems in geometric structures (see
Gallo et al.~\cite{Gallo/Kapovich/Marden:2000}
where the authors
refer to these geometries as \emph{branched}
rather than \emph{generalized}).
The \Thullen and Hartogs extension problems for maps to complex homogeneous
spaces (which can be studied using invariants derived
via Cartan's method of the moving frame)
are largely untouched, besides the work of Ivashkovich~\cite{Ivashkovich:2008}.
The obvious meromorphic extension problems for meromorphic
geometric structures are open.

\subsection{Curvature of Hermitian metrics and extension problems}

It is easy to generalize many of the theorems of this paper
to utilize integral curvature bounds.

\begin{Theorem}[Shevchishin~\cite{Shevchishin:1996}]
A holomorphic vector bundle on an $n$-manifold
extends across a subset of \cct just when it admits a Hermitian
metric with $L^n$-bounded curvature.
\end{Theorem}

This result in particular characterizes the nondegenerate
plane f\/ields that extend across subsets of \cct,
although the characterization is
not easy to apply to examples, since we have to pick
a metric for which we can see how to bound the integral.

\subsection*{Acknowledgments}

This material is based upon works supported by the Science Foundation
Ireland under Grant No. MATF634. The author is grateful to the
editors, and to Sorin Dumitrescu, Sergei Ivashkovich and 
anonymous reviewers
for their assistance in improving this paper.

\addcontentsline{toc}{section}{References}
\LastPageEnding

\end{document}